\documentclass[12pt,a4paper]{amsart}

\usepackage[hmargin=2cm,vmargin=2cm]{geometry}

\usepackage{tikz}

\usepackage{amsmath,amsfonts,amssymb}

\def\E{\mathrm{E}}

\def\L{\mathrm{L}}
\def\V{\mathrm{V}}

\newcommand{\AAA}{{\mathsf{A}}}
\newcommand{\BBB}{{\mathsf{B}}}
\newcommand{\CCC}{{\mathsf{H}}}
\newcommand{\B}{\CCC}
\newcommand{\tP}{{\mathcal{P}}}
\newcommand{\cR}{{\mathcal{R}}}
\newcommand{\cS}{{\mathcal{S}}}

\newcommand{\Z}{\mathbb Z}
\newcommand{\Q}{\mathbb Q}
\newcommand{\mbR}{\mathbb R}
\newcommand{\mbC}{\mathbb C}
\def\CP1{\mathbb{CP}^1}
\newcommand{\com}{{\mathbb C}}
\newcommand{\MM}{{\mathcal{M}}}

\newcommand{\mcM}{\mathcal M}
\newcommand{\oM}{\overline{\mathcal M}}

\newcommand{\cL}{\mathcal L}

\newcommand{\mcH}{\mathcal H}
\newcommand{\cF}{\mathcal F}
\newcommand{\Mat}{\mathop{\mathrm{Mat}}\nolimits}
\renewcommand{\Im}{\mathop{\mathrm{Im}}\nolimits}
\renewcommand{\Re}{\mathop{\mathrm{Re}}\nolimits}
\newcommand{\tr}{\mathop{\mathrm{tr}}\nolimits}
\newcommand{\Ai}{\mathop{\mathrm{Ai}}\nolimits}
\newcommand{\ev}{\mathop{\mathrm{ev}}\nolimits}
\newcommand{\End}{\mathop{\mathrm{End}}\nolimits}

\newcommand{\tF}{\widetilde F}

\newcommand{\<}{\left <}
\renewcommand{\>}{\right >}

\newcommand{\vac}{\left|0\right>}

\def\d{\partial}

\def\s{\mathfrak{s}}

\newcommand{\half}{\frac{1}{2}}

\def\Coef{{\rm Coef}}
\def\Gr{{\rm Gr}}

\newtheorem{theorem}{Theorem}

\newtheorem{question}{Question}

\theoremstyle{definition}

\begin{document}

\title[]{The hypergeometric functions of the
Faber-Zagier and Pixton relations}
\author {A. Buryak}
\address {ETH Z\"urich, Department of Mathematics}
\email {alexander.buryak@math.ethz.ch}
\author {F. Janda}
\address {ETH Z\"urich, Department of Mathematics}
\email {janda@math.ethz.ch}
\author{R. Pandharipande}
\address{ETH Z\"{u}rich, Department of Mathematics}
\email {rahul@math.ethz.ch}
\date{July 2015}

\maketitle

\begin{abstract} The relations in the tautological
ring of the moduli space $\MM_{g}$ of nonsingular curves 
conjectured by Faber-Zagier in 2000 and extended 
to the moduli space $\overline{\MM}_{g,n}$ of stable curves
by Pixton in 2012 are based upon two hypergeometric
series $\AAA$ and $\BBB$. The question of the
geometric origins of these series has been solved in at least
two ways (via the Frobenius structures associated to 3-spin curves
and to $\mathbb{P}^1$). The series $\AAA$ and $\BBB$ 
also appear in the study of descendent integration
on the moduli spaces of open and
closed curves. We survey here 
the various occurrences of $\AAA$ and $\BBB$ starting from their
appearance in 
the asymptotic expansion of the Airy function (calculated by Stokes
in the $19^{th}$ century). Several open questions are proposed.
\end {abstract}

\setcounter{tocdepth}{1} 
\tableofcontents

\setcounter{section}{-1}

\section{Introduction}

\subsection{Tautological classes}
For $g\geq 2$, let $\MM_g$ be the moduli space of nonsingular, projective,
genus $g$ curves over $\com$, and let
\begin{equation*}
\pi: \mathcal{C}_g \rightarrow \mathcal{M}_g 
\end{equation*}
be the universal curve. 
The cotangent line class is defined via the
line bundle  $\omega_\pi$ of relative differentials of the 
morphism $\pi$,
$$\psi = c_1(\omega_\pi)\in A^1(\mathcal{C}_g,\mathbb{Q})\ .$$
The $\kappa$ classes are defined by push-forward,
$$\kappa_r = \pi_*( \psi^{r+1}) \in A^{r}(\MM_g,\mathbb{Q})\ .$$
The {\em tautological ring} in algebraic cycles,
$$R^*(\MM_g) \subset A^*(\MM_g, \mathbb{Q})\, ,$$
is the $\mathbb{Q}$-subalgebra generated by all of the
$\kappa$ classes. 
Since 
$$\kappa_{0}= 2g-2 \in \mathbb{Q}$$
is a multiple of the fundamental class, we need not take 
$\kappa_0$ as a generator.
There is a canonical quotient
$$\mathbb{Q}[\kappa_1,\kappa_2, \kappa_3, \ldots] 
\stackrel{q}{\longrightarrow} R^*(\MM_g) \longrightarrow 0\ .$$
The ideal of {\em tautological relations} among the
$\kappa$ classes is the kernel of $q$.

\subsection{Relations}
Faber and  Zagier conjectured in 2000 a remarkable set
of relations among the $\kappa$ classes in $R^*(\MM_g)$
which were first proven to hold in \cite{PP13}.

To write the Faber-Zagier relations, we will require
the following notation.
Let the variable set
\vspace{-5pt}
$$\mathbf{p} = \{\ p_1,p_3,p_4,p_6,p_7,p_9,p_{10}, \ldots\ \}$$
be indexed by positive integers {\em not} congruent
to $2$ modulo $3$.
Define the series
\begin{eqnarray*}
\Psi(t,\mathbf{p})& = &
\left(1+tp_3+t^2p_6+t^3p_9+\ldots\right) \sum_{i=0}^\infty \frac{(6i)!}{(3i)!(2i)!} t^i
\\ & & +\ \left(p_1+tp_4+t^2p_7+\ldots\right ) 
\sum_{i=0}^\infty \frac{(6i)!}{(3i)!(2i)!} \frac{6i+1}{6i-1} t^i \ .
\end{eqnarray*}
Since $\Psi$ has constant term 1, we may take the logarithm.
Define the constants $C_r^{\text{\tiny{{\sf FZ}}}}(\sigma)$ by the formula
$$\log(\Psi)= 
\sum_{\sigma}
\sum_{r=0}^\infty C_r^{\text{\tiny{{\sf FZ}}}}(\sigma)\ t^r 
\mathbf{p}^\sigma
\ . $$
The above sum is over all partitions $\sigma$ of size 
$|\sigma|$ 
which avoid 
 parts congruent to 2 modulo 3. The empty partition is included
in the sum.
To the partition  $\sigma=1^{n_1}3^{n_3}4^{n_4} \cdots$, we associate
the monomial
$\mathbf{p}^\sigma= p_1^{n_1}p_3^{n_3}p_4^{n_4}\cdots$.
Let 
$$\gamma^{\text{\tiny{{\sf FZ}}}}
= 
\sum_{\sigma}
 \sum_{r=0}^\infty C_r^{\text{\tiny{{\sf FZ}}}}(\sigma)
\ \kappa_r t^r 
\mathbf{p}^\sigma
\ .
$$
For a series $\Theta\in \mathbb{Q}[\kappa][[t,\mathbf{p}]]$ in the variables $
\kappa_i$, $t$, and $p_j$, let
$[\Theta]_{t^r \mathbf{p}^\sigma}$ denote the
 coefficient of $t^r\mathbf{p}^\sigma$
(which is a polynomial in the $\kappa_i$).

\begin{theorem}[Pandharipande-Pixton] \label{dddd}
{ In $R^r(\MM_g)$, the Faber-Zagier relation
$$
\big[ \exp(-\gamma^{\text{\tiny{{\sf FZ}}}}) \big]_{t^r \mathbf{p}^\sigma}  = 0$$
holds when
$g-1+|\sigma|< 3r$ and
$g\equiv r+|\sigma|+1 \mod 2$.}
\end{theorem}

As a corollary of the proof \cite{PP13} of Theorem \ref{dddd}, 
a stronger boundary result was obtained.
If $g-1+|\sigma|< 3r$ and
$g\equiv r+|\sigma|+1 \mod 2$, then
\begin{equation}
\big[ \exp(-\gamma^{\text{\tiny{{\sf FZ}}}}) \big]_{t^r \mathbf{p}^\sigma}  \in R^*(\partial\overline{\mathcal{M}}_g)\ .
\end{equation}
Not only is the Faber-Zagier relation 0 in $R^*(\mathcal{M}_g)$, but the
relation is equal to a tautological class supported on the boundary of
the moduli space $\overline{\mathcal{M}}_g$. 

A precise conjecture for 
the boundary terms (and much more) has been proposed 
by Pixton in \cite{Pix12}.
We review the complete form of Pixton's relations in Appendix A,
see also \cite{PPZ13,Pix12}. Pixton has conjectured that his relations
provide a {\em complete} set of tautological relations in the Chow
rings of the moduli spaces 
$\overline{\MM}_{g,n}$. Since  Pixton's relations restrict to
the Faber-Zagier relations, the hypergeometric series 
$$\AAA(z)=\sum_{i=0}^\infty \frac{(6i)!}{(3i)!(2i)!} 
\left(\frac{z}{288}\right)^i\, , \ \ \ \
\BBB(z)=\sum_{i=0}^\infty \frac{(6i)!}{(3i)!(2i)!} \frac{6i+1}{6i-1} 
\left(\frac{z}{288}\right)^i\, , \ \ \  z=288t $$
also occur in the formula of Pixton.
In fact, just as above for Faber and Zagier, 
the series $\AAA$ and $\BBB$ are the {\em only} non-formal
inputs for Pixton.

\subsection{Differential equations}
The series $\AAA$ and $\BBB$ are easily related
via the following differential equation:
$$3 z^2\ \frac{d\AAA}{dz} +
\left(\frac{z}{2}-1\right)\AAA =\BBB\, \ .$$
Hence, we often view $\AAA$ as the more fundamental function.
The main hypergeometric differential equation satisfied
by $\AAA$ is:
$$3 z^2 \frac{d^2\AAA}{dz^2}  + (6z-2) \frac{d\AAA}{dz} + \frac{5}{12} \AAA = 0 \ .$$

\subsection{Origins of $\AAA$ and $\BBB$}
In order to prove the Faber-Zagier and Pixton 
relations, geometric sources for the
series $\AAA$ and $\BBB$ were found. At present,  two successful
approaches are known: via
the Frobenius geometries of 
$3$-spin 
curves \cite{PPZ13} and of
$\mathbb{P}^1$ \cite{Jan13,Jan14,PP13}. The two approaches 
lead to two different geometric origins for $\AAA$ and $\BBB$.

More recently, occurances
of $\AAA$ have been noticed \cite{Bur14b}
in the generating series of descendent integrals over
the moduli spaces of open Riemann surfaces \cite{PST}. Remarkably,
the series $\AAA$ can already be seen in the 
asymptotic expansion of the Airy function related to the Witten-Kontsevich theory
of descendent integration over $\overline{\MM}_{g,n}$.
Our goal here is to survey these various 
appearances of the series $\AAA$ and $\BBB$.

The  occurances
connected to descendent integration have not (yet)
played a role in proofs of the Faber-Zagier and Pixton relations.
Perhaps the reverse is more likely: the relations could be
used to
constrain descendent integration. For integration against
the product of the top two Chern classes of the Hodge bundle,
$$\lambda_g\lambda_{g-1}\in A^{2g-1}(\overline{\MM}_{g,n},\mathbb{Q})\, ,$$
corresponding to the geometry of the moduli space of nonsingular curves,
a subset of Pixton's relations have been shown in \cite{PPZ15,PX} to imply the
$\lambda_g\lambda_{g-1}$ descendent formula \cite{Fa,GP}. 

\subsection{Acknowledgements}
We thank C. Faber, A. Pixton, S. Shadrin, J. Solomon, R. Tessler,
and D. Zvonkine 
for discussions related to the tautological ring,
descendent integration, and the hypergeometric series $\AAA$ and $\BBB$.

A.B. was supported by the grants ERC-2012-AdG-320368-MCSK, RFFI 13-01-00075,
and NSh-4850.2012.1.
F.J. was supported by the grant SNF-200021-143274.
R.P. was partially supported by SNF-200021-143274, ERC-2012-AdG-320368-MCSK, 
SwissMap, and the Einstein Stiftung.
The paper was completed while F.J. and R.P. were visiting 
Humboldt University in Berlin.

\section{Asymptotic expansion of the Airy function} \label{sec:asymp}
Define the closely related functions $\mathcal{A}$ and $\mathcal{B}$
by
$$\mathcal{A}(x)= \AAA(-x^3)\ , \ \ \ \mathcal{B}(x)= \BBB(-x^3)\, .$$
We see
$$\mathcal{A}(x) =\sum_{j= 0}^\infty a_j x^{3j}= 1-\frac{5}{24}x^3 +\ldots\, , \ \ \ \ 
-\mathcal{B}(x)=\sum_{j= 0}^\infty b_j x^{3j} = 1+\frac{7}{24}x^3+\ldots \, .$$
where
the coefficients $a_j$ and $b_j$ are
$$a_j=(-1)^j\frac{(6j)!}{288^j(2j)!(3j)!}\, , \ \ \ \
b_j=-\frac{6j+1}{6j-1}a_j\,  \ \ \ \ j\geq 0\, .$$

The series $\mathcal{A}(x)$ is related to the classical Airy function in the following way. The Airy function $\Ai(x)$ is defined by
\begin{equation}
  \label{eq:defairy}
  \Ai(x)=\int_0^{\infty}\cos\left(\frac{t^3}{3}+xt\right)dt\, ,\quad x\in\mbR\, .
\end{equation}
It is the unique (up to a scalar factor) bounded real solution of the
Airy differential equation
\begin{equation}\label{tqq1}
y''=x y\, .
\end{equation}
The Airy function $\Ai(x)$ has the following asymptotic expansion for $x\to\infty$, 
\begin{equation}\label{tt33}
\Ai(x)\, \asymp_{x\to\infty}\, \frac{\sqrt\pi}{2}x^{-\frac{1}{4}}e^{-\frac{2}{3}x^{\frac{3}{2}}}
\mathcal{A}\left(2^{-\frac{1}{3}}x^{-\half}\right),
\end{equation}
see \cite[pages 22-23]{EMOT53}.

We review here the short derivation of asymptotic expansion \eqref{tt33}.
We may write the oscillatory integral defining the Airy function as
\begin{equation} \label{vv22}
 \Ai(x) = \int_0^\infty \cos\left(\frac{t^3}3 + xt\right) \mathrm dt
  = \frac 12 \int_{-\infty}^\infty e^{i\left(\frac{t^3}3 + xt\right)} \mathrm dt\, .
\end{equation}
Viewing \eqref{vv22} as a complex line integral, we move the
integration contour from $\{\Im(t) = 0\}$ to $\{\Im(t) = i\sqrt{x}\}$
by shifting the integration variable $t$
by $i\sqrt{x}$.
There are no poles in the region 
$$0\le\mathop{\mathrm{Im}}t\le\sqrt{x}\,.$$
 Let us check that the integrals of the function~$\frac{1}{2}e^{i\left(\frac{t^3}{3}+xt\right)}$ over the arcs~$\Gamma_1$ and~$\Gamma_2$ (see Figure ~\ref{fig:shift}) go to zero, when~$R$ goes to infinity. If~$R$ is big enough, then $\alpha\le\frac{\pi}{6}$ and we have 
\begin{eqnarray*}
\left|\frac{1}{2}\int_{\Gamma_2}e^{i\left(\frac{t^3}{3}+xt\right)}dt\right|&=&\left|\frac{1}{2}\int_{\Gamma_1}e^{i\left(\frac{t^3}{3}+xt\right)}dt\right|\\
&\le & \frac{R}{2}\int_0^{\alpha}\left|e^{i\left(\frac{R^3}{3}e^{3i\phi}+x R e^{i\phi}\right)}\right|d\phi\\
&= &\frac{R}{2}\int_0^{\alpha}e^{-\frac{R^3}{3}\sin 3\phi-x R \sin\phi}d\phi\\
&\stackrel{\substack{\text{by Jordan's}\\\text{inequality}}}{\le} &\frac{R}{2}\int_0^{\alpha}e^{-\frac{2}{\pi}(R^3+x R)\phi}d\phi\\
&=&\frac{\pi}{4(R^2+x)}\left(1-e^{-\frac{2}{\pi}(R^3+xR)\alpha}\right)\xrightarrow[R\to\infty]{}0\, .
\end{eqnarray*}
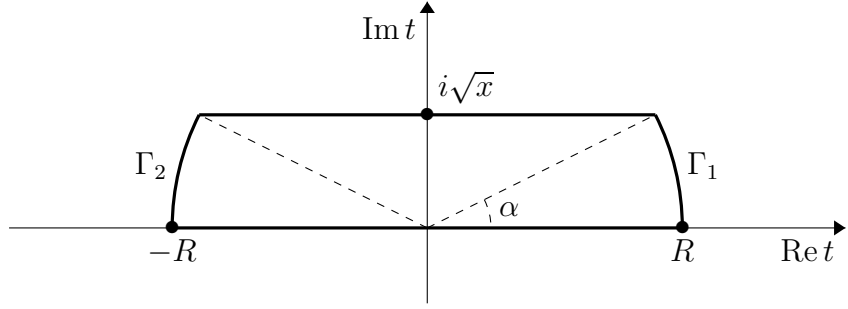
\begin{figure}[t]

\begin{tikzpicture}[scale=0.5]

\draw (-11,0) -- (11,0);
\fill (10.7,0.2)--(11,0)--(10.7,-0.2)--cycle;
\draw (0,-2) -- (0,6);
\fill (-0.2,5.7)--(0,6)--(0.2,5.7)--cycle;

\draw [very thick] (6,3) arc (26.57:0:6.71);
\draw [very thick] (-6,3) arc (153.43:180:6.71);
\draw [very thick] (-6.71,0) -- (6.71,0);
\draw [very thick] (-6,3) -- (6,3);

\draw [dashed] (0,0) -- (6,3);
\draw [dashed] (0,0) -- (-6,3);
\draw [dashed] (1.5,0.75) arc (26.57:0:6.71/4);

\coordinate [label=-90: $\mathop{\mathrm{Re}}t$] (B) at (10,0);
\coordinate [label=180: $\mathop{\mathrm{Im}}t$] (B) at (0,5.3);
\coordinate [label=45: $i\sqrt{x}$] (B) at (0,3);
\coordinate [label=center: \textbullet] (B) at (0,3);
\coordinate [label=center: \textbullet] (B) at (6.71,0);
\coordinate [label=center: \textbullet] (B) at (-6.71,0);
\coordinate [label=-90: $R$] (B) at (6.71,0);
\coordinate [label=-90: $-R$] (B) at (-6.71,0);
\coordinate [label=0: $\alpha$] (B) at (1.6,0.5);
\coordinate [label=0: $\Gamma_1$] (B) at (6.55,1.65);
\coordinate [label=180: $\Gamma_2$] (B) at (-6.55,1.65);
\end{tikzpicture}

\caption{Shift of the integration contour}
\label{fig:shift}
\end{figure}
Therefore, we obtain
\begin{equation*}
  \Ai(x) = \frac 12 \int\limits_{-\infty}^\infty e^{i \frac{t^3}3 - \sqrt{x} t^2 - \frac 23 x^{3/2}} \mathrm dt\, .
\end{equation*}
Scaling $t$ by $(2\sqrt{x})^{1/2}$ makes the integrand a deformed
Gaussian integral:
\begin{equation*}
  \Ai(x) = e^{-\frac 23 x^{3/2} }\, \frac{x^{-\frac{1}{4}}}{2\sqrt{2}}\, \int\limits_{-\infty}^\infty e^{-\frac{t^2}2} e^{i x^{-3/4} \frac{t^3}{6\sqrt{2}}} \, \mathrm dt\, .
\end{equation*}
By dominated convergence, we see 
\begin{equation*}
  e^{-\frac 23 x^{3/2} } \frac{x^{-\frac{1}{4}}}{2\sqrt{2}} \sum_{j = 0}^\infty \int\limits_{-\infty}^\infty e^{-\frac{t^2}2} \frac 1{j!} \left(\frac{i t^3}{6\sqrt{2} x^{3/4}}\right)^j \mathrm dt
  = e^{-\frac 23 x^{3/2} } \frac{x^{-\frac{1}{4}}}{2\sqrt{2}} \sum_{j = 0}^\infty \frac{x^{-3j/4}}{(-6\sqrt{2}i)^j j!} \int\limits_{-\infty}^\infty e^{-\frac{t^2}2} t^{3j} \mathrm dt
\end{equation*}
is an asymptotic expansion for $\Ai(x)$ for $x \to \infty$. Using 
\begin{equation*}
  \int\limits_{-\infty}^\infty e^{-\frac{t^2}2} t^j \mathrm dt =
  \begin{cases}
    \sqrt{2\pi}\, (j - 1)!!\, , & \text{if $j$ is even,} \\
    0\, , & \text{if $j$ is odd,}
  \end{cases}
\end{equation*}
we arrive at
\begin{eqnarray*}
  \Ai(x) &\asymp_{x\rightarrow \infty}& \frac{\sqrt{\pi}}2 e^{-\frac 23 x^{3/2} } 
x^{-\frac{1}{4}} \sum_{j = 0}^\infty \frac{(6j - 1)!!}{(6\sqrt{2}i)^{2j} (2j)!} x^{-3j/2}
 \\ 
& & \ \ \ \  =  \frac{\sqrt{\pi}}2 e^{-\frac 23 x^{3/2} } x^{-\frac{1}{4}} \sum_{j = 0}^\infty \frac{(6j)!}{(2j)!(3j)!} \left(-\frac{x^{-3/2}}{576}\right)^j \, ,
\end{eqnarray*}
which is exactly \eqref{tt33}.

Similarly, for the derivative of $\Ai(x)$, 
\begin{equation*}
 \Ai'(x) =  \frac 12 \int\limits_{-\infty}^\infty e^{i\left(\frac{t^3}3 + xt\right)} it\,
 \mathrm dt\, ,
\end{equation*}
we obtain the asymptotic expansion
\begin{equation*}
\Ai'(x) 
\, \asymp_{x\rightarrow \infty}\, \frac{\sqrt{\pi}}2 e^{-\frac 23 x^{3/2} } x^{\frac{1}{4}} \sum_{j = 0}^\infty \frac{(6j)!}{(2j)!(3j)!} \frac{6j + 1}{6j - 1} \left(-\frac{x^{-3/2}}{576}\right)^j\, , 
\end{equation*}
or equivalently,
\begin{equation*}
\Ai'(x)\, \asymp_{x\to\infty}\, \frac{\sqrt\pi}{2}x^{\frac{1}{4}}e^{-\frac{2}{3}x^{\frac{3}{2}}}
\mathcal{B}\left(2^{-\frac{1}{3}}x^{-\half}\right)\, .
\end{equation*}


\section{Moduli of stable curves and the
 infinite Grassmanian}\label{section:closed curves}
\subsection{Overview}
The series $\mathcal{A}$ appears in the intersection theory of the 
moduli space $\overline{\MM}_{g,n}$ of stable curves.
We start by reviewing Witten's conjecture 
governing descendent
integration in 
Section~\ref{subsection:Witten's conjecture}.
Kontsevich's formula for the generating series $F^c$ of the descendents 
is expressed in terms of integrals over spaces of Hermitian matrices. 
Certain specializations of the  descendent partition function  $\exp(F^c)$ 
may be
expressed as ratios of simple determinants. 
As an almost immediate consequence, 
the most  basic of these specialization of  $\exp(F^c)$ 
coincides with the hypergeometric series $\mathcal{A}$. 

In Section~\ref{subsection:tau-functions}, we review 
a well-known construction associating a tau-function of the KP hierarchy 
to any point of the infinite dimensional Grassmanian. 
We then present a result of Kac and Schwarz which explicitly 
describes the point in the Grassmanian corresponding to the 
partition function $\exp(F^c)$. 
The series $\mathcal{A}$ and $\mathcal{B}$ emerge here and play a prominent role.  

\subsection{Witten's conjecture}\label{subsection:Witten's conjecture}

Let $\oM_{g,n}$ be the moduli space of genus $g$ stable curves over
$\mathbb{C}$ with $n$ 
marked points. The first Chern class of the
contangent line at the $i^{th}$ marking is denoted by
$$\psi_i\in A^1(\oM_{g,n},\mathbb{Q})\, .$$ 
We define the {\em descendent integrals} by
\begin{equation}
\<\tau_{k_1}\tau_{k_2}\ldots\tau_{k_n}\>^c_g=
\int_{\oM_{g,n}}\psi_1^{k_1}\psi_2^{k_2}\ldots\psi_n^{k_n}\, .
\label{frss}
\end{equation}
The bracket \eqref{frss} vanishes unless the
dimension constraint
$$3g-3+n = \sum_{i=1}^n k_i$$
is satisfied.
The superscript $c$ here indicates
 integration over the moduli of compact Riemann surfaces. 
Integration over the moduli of open Riemann surfaces
will be considered in Section~\ref{section:open curves}. 

We introduce variables $\{t_i\}_{i\geq 0}$ and define the generating series 
$$
F^c(t_0,t_1,\ldots)=\sum_{\substack{g\geq 0,\,  n\geq 1\\ \vspace{-2pt} \\ 2g-2+n>0}}\frac{1}{n!}\sum_{k_1,\ldots,k_n\ge 0}\<\tau_{k_1}\tau_{k_2}\ldots\tau_{k_n}\>^c_gt_{k_1}t_{k_2}
\ldots t_{k_n}\, .
$$
Witten \cite{Wit91} conjectured that the 
partition function $\exp(F^c)$ is a tau-function of the KdV hierarchy. 
In particular, 
$$u=\frac{\d^2 F^c}{\d t_0^2}$$ is a solution of the KdV hierarchy. 
The first equations of the hierarchy are
\begin{align*}
u_{t_1}&=uu_x+\frac{1}{12}u_{xxx},\\ \notag
u_{t_2}&=\frac{1}{2}u^2u_x+\frac{1}{12}(2u_xu_{xx}+uu_{xxx})
+\frac{1}{240}u_{xxxxx}\, ,\\ \notag
& \vdots
\end{align*}
where we have identified $x$ here with $t_0$. 
Witten's conjecture was proven by Kontsevich \cite{Kon92}.
See \cite{KL07,Mir,OP1} for other proofs. 


\subsection{Kontsevich's matrix integral}\label{subsection:matrix integral}

\subsubsection{Matrix integrals}\label{subsubsection:matrix integration}
Kontsevich \cite{Kon92} proposed a representation of the partition function $\exp(F^c)$ in terms of
integrals over spaces of Hermitian matrices. To connect $\exp(F^c)$ to
the hypergeometric series $\AAA$, we will use 
Kontsevich's matrix model.

Let $\mcH_N$ denote the $N^2$-dimensional real vector space of Hermitian $N\times N$ matrices,
$$
\mcH_N=\Big\{\ H=(h_{i,j})\in\Mat_{N,N}(\mbC)\ \Big| \ h_{i,j}=\overline{h}_{j,i}\ \Big\}.
$$
We introduce coordinates $\{x_{i,i}\}_{1\le i\le N}$ and $\{x_{i,j},y_{i,j}\}_{1\le i<j\le N}$ on $\mcH_N$ by
\begin{align*}
&x_{i,i}=h_{i,i},\\
&x_{i,j}=\Re(h_{i,j})\quad\text{and}\quad y_{i,j}=\Im(h_{i,j}),\quad i<j. 
\end{align*}
A volume form on $\mcH_N$ is defined by
$$
dv(H)=\prod_{i} d x_{i,i}\,\prod_{i<j}d x_{i,j}d y_{i,j}\, .
$$

Let $\Lambda$ be a diagonal $N\times N$ matrix with positive real entries $\Lambda_1,\ldots,\Lambda_N$ along
 the diagonal. Define a Gaussian measure on the space $\mcH_N$ by
$$d\mu_\Lambda(H)=c_{\Lambda,N}\, e^{-\half\tr(H^2\Lambda)}\, dv(H),$$
where the normalization 
$$c_{\Lambda,N}=(2\pi)^{-N^2/2}\, \prod_{i}\Lambda_i^{1/2}\, \prod_{i<j}(\Lambda_i+\Lambda_j)$$
is determined by the constraint
$$\int_{\mcH_N}d\mu_{\Lambda,N}(H)= 1 \, .$$ 


Since $\left|e^{\frac{i}{6}\tr(H^3)}\right|\le 1$, for $H\in\mcH_N$, the integral
\begin{equation}\label{gvvt}
\int_{\mcH_N} e^{\frac{i}{6}\tr(H^3)}d\mu_{\Lambda,N}(H)
\end{equation}
is convergent and defines a function of $(\Lambda_1,\ldots,\Lambda_N)\in\mbR_{>0}^N$. 
When $\Lambda_j^{-1}\to 0$, the function \eqref{gvvt} admits an asymptotic expansion given by
\begin{gather}\label{eq:Kontsevich's formula}
\int_{\mcH_N} e^{\frac{i}{6}\tr(H^3)}d\mu_{\Lambda,N}(H)\ \asymp_{\Lambda_j^{-1}\to 0}\ \left.\exp(F^c)\right|_{t_i=-(2i-1)!!\tr(\Lambda^{-2i-1})}\, ,
\end{gather}
The above expansion \eqref{eq:Kontsevich's formula} is called 
Kontsevich's formula \cite{Kon92}.

\subsubsection{Determinantal formulas}
Via an averaging procedure over the unitary group applied 
to the left side of \eqref{eq:Kontsevich's formula}, the following formula
can be obtained:
\begin{gather}\label{eq:determinantal formula}
\left.\exp(F^c)\right|_{t_i=-(2i-1)!!\tr(\Lambda^{-2i-1})}=\frac{\left|(D_j^{i-1}\mathcal A(\Lambda^{-1}_j))_{1\le i,j\le N}\right|}{\prod_{1\le i<j\le N}(\Lambda_j-\Lambda_i)}\, ,
\end{gather}
where $D_j$ is the differential operator 
$$
D_j=-\frac{1}{\Lambda_j}\frac{d}{d\Lambda_j}+\Lambda_j+\frac{1}{2\Lambda_j^2}\, .
$$
We recommend \cite[Section 2.2]{IZ92} for a quick derivation of
\eqref{eq:determinantal formula}. 
In Section \ref{subsubsection:special case}, we provide a direct proof for $N=1$.

\subsubsection{Case $N=1$}\label{subsubsection:special case}

For $N=1$, Kontsevich's formula~\eqref{eq:Kontsevich's formula} yields
$$
\frac{\sqrt{\lambda}}{\sqrt{2\pi}}\int_\mbR e^{\frac{i}{6}t^3-\half t^2\lambda}dt\ \asymp_{\lambda^{-1}\to 0}\ 
\left.\exp(F^c)\right|_{t_i=-(2i-1)!!\lambda^{-2i-1}}.
$$
Clearly, the left-hand side is equal to $\frac{1}{\sqrt{2\pi}}\int_\mbR e^{\frac{i}{6}\lambda^{-\frac{3}{2}}t^3}e^{-\half t^2}dt$. We compute:
\begin{eqnarray*}
\frac{1}{\sqrt{2\pi}}\int_\mbR e^{\frac{i}{6}\lambda^{-\frac{3}{2}}t^3}e^{-\half t^2}dt
& \asymp_{\lambda^{-1}\to 0}& \sum_{j=0}^{\infty}\left(\frac{1}{\sqrt{2\pi}}\int_{\mbR}\frac{i^j}{6^j j!}t^{3j}e^{-\half t^2}dt\right)\lambda^{-\frac{3}{2}j}\\
& & =\sum_{j=0}^{\infty}(-1)^j\frac{(6j-1)!!}{36^j(2j)!}\lambda^{-3j}\\ & &  =\mathcal{A}(\lambda^{-1}).
\end{eqnarray*}
We have arrived at a direct connection between descendent integration and the series $\mathcal{A}(\lambda^{-1})$, 
\begin{equation}
\label{qqq}
\left.\exp(F^c)\right|_{t_i=-(2i-1)!!\lambda^{-2i-1}}=\mathcal{A}(\lambda^{-1}).
\end{equation}

Pixton's relations constrain tautological classes on 
$\overline{\mathcal{M}}_{g,n}$ and hence also descendent integrals.
In fact,
Pixton's relations are expected to uniquely determine the
descendent theory, but the implication is not yet proven.
A simpler question, since both involve the
hypergeometric series $\mathcal{A}$, is the following.

\begin{question}
  Can the specialization of the partition fuction \eqref{qqq} be
  derived from Pixton's relations?
\end{question}


\pagebreak

\subsection{The Infinite Grassmanian and tau-functions of the KP hierarchy}\label{subsection:tau-functions}



\subsubsection{Brief introduction to the KP hierarchy}\label{subsubsection:brief introduction}

A {\em pseudo-differential} operator $A$ is a Laurent series 
$$
A=\sum_{n=-\infty}^m a_n(T)\, \d_x^n\, ,
$$
where $m\in \mathbb{Z}$ and the coefficients $a_n(T)$ are formal power series in the
variables 
$\{T_i\}_{i\geq 1}$,
$$a_n(T) \in \mathbb{C}[[T_1,T_2,T_3, \ldots]]\, .$$
We identify the variable $x$ with $T_1$. 
The non-negative and negative degree parts of 
the pseudo-differential operator $A$ are defined by
\begin{gather*}
A_+ =\sum_{n=0}^m a_n\d_x^n \quad\text{and}\quad A_-=A-A_+ \, .
\end{gather*}

The product of pseudo-differential operators is defined by the following commutation rule:
\begin{gather*}
\d_x^k\circ f=\sum_{l=0}^\infty\frac{k(k-1)\ldots(k-l+1)}{l!}\frac{\d^lf}{\d x^l}\d_x^{k-l} \, ,
\end{gather*} 
where $k\in\Z$ and $f\in\mbC[[T_1,T_2,T_3,\ldots]]$.

Consider the pseudo-differential operator
$$
L=\d_x+\sum_{i\ge 1}w_i\d_x^{-i}\, .
$$
The {\em KP hierarchy} is the following system of partial differential equations for the power series~$w_i$:
\begin{gather}\label{eq:KP hierarchy}
\frac{\d L}{\d T_n}=\left[\left(L^n\right)_+,L\right],\quad n=1,2,3,\ldots.
\end{gather}
For $n=1$, the equation is equivalent to
 $$\frac{\d w_i}{\d T_1}=\frac{\d w_i}{\d x}\, , \ \ \forall i\geq 1\, ,$$ 
compatible with our identification of $x$ with $T_1$.

Suppose an operator $L$ satisfies the system~\eqref{eq:KP hierarchy}. Then there exists a pseudo-differential operator $P$ of the form
\begin{gather}\label{dressing operator}
P=1+\sum_{n\ge 1}p_n(T)\, \d_x^{-n},
\end{gather}
satisfying $L=P\circ\d_x\circ P^{-1}$ and
\begin{align}
&\frac{\d P}{\d T_n}=-\left(L^n\right)_-\circ P,\quad\ \ n=1,2,3,\ldots\ .\label{eq:Sato-Wilson equations}
\end{align}
The operator $P$ is the {\em dressing} operator and 
\eqref{eq:Sato-Wilson equations} are the {Sato-Wilson equations}. The Laurent series 
$$
\widehat P(T;z)=1+\sum_{n\ge 1}p_n(T)\,z^{-n}
$$
is the {\em symbol} of the dressing operator $P$.

We can now introduce the notion of a tau-function. 
Denote by $G_z$ the {\em shift operator} which acts on 
a power series $f\in\mbC[[T_1,T_2,T_3,\ldots]]$ as follows:
\begin{gather*}
G_z(f)\Big(T_1,T_2,T_3,\ldots\Big)=f\Big(T_1-\frac{1}{z},\, 
T_2-\frac{1}{2z^2},\, T_3-\frac{1}{3z^3},\ldots\Big)\, .
\end{gather*}
Let $P=1+\sum_{n\ge 1}p_n(T)\,\d_x^{-n}$ be the dressing operator of some operator $L$ satisfying the KP hierarchy~\eqref{eq:KP hierarchy}. 
Then there exists a series $\tau\in\mbC[[T_1,T_2,T_3,\ldots]]$ 
with constant term $\left.\tau\right|_{\{T_i=0\}}=1$ for which
$$
\widehat P=\frac{G_z(\tau)}{\tau}\, .
$$
The series $\tau$ is a {\em tau-function} of the KP hierarchy.

The KdV hierarchy is a certain reduction of the KP hierarchy. 
We do not discuss the details here, but only state the following
property: a tau-function of the KdV hierarchy is a tau-function of the KP hierarchy which is independent of the variables $\{T_{2i}\}_{i\geq 1}$.

The precise form of Witten's conjecture may now be formulated: 
$$
\left.\exp(F^c)\right|_{t_i=(2i+1)!!T_{2i+1}}
$$
is a tau-function of the KdV hierarchy.

\subsubsection{The Infinite Grassmanian and the Fock space}\label{subsubsection:Grassmanian}

Consider the space of Laurent series $\mbC[z^{-1},z]]$. There is a natural projection 
$$p_-\colon\mbC[z^{-1},z]]\to\mbC[z^{-1}]\, .$$
 We denote by $\Gr^0\mbC[z^{-1},z]]$ the set of all vector subspaces $H\subset\mbC[z^{-1},z]]$ for which the projection $$p_-\colon H\to\mbC[z^{-1}]$$
 is an isomorphism. 

The {\em Fock space} $\cF$ is the vector space of (possibly infinite)
linear combinations of wedge products of the form
\begin{gather}\label{Fock basis}
z^{a_1}\wedge z^{a_2}\wedge z^{a_3} \wedge \ldots,
\end{gather}
where $a_1,a_2,a_3,\ldots$ is a decreasing sequence of integers 
for which there exists an integer $c$ (called the charge) satisfying
 $$a_i=-i+c$$ for all $i$ sufficiently large. 
Denote by $\cF^{[c]}\subset\cF$ the subspace consisting of vectors with charge $c$.

The {\em vacuum} vector 
$$
\vac =z^{-1}\wedge z^{-2}\wedge z^{-3}\wedge z^{-4} \wedge z^{-5}\ldots\ \ \in\cF^{[0]}
$$
plays a special role. 
For any vector $v\in\cF$, we  denote by $\left<0\,|\,v\right>$ the coefficient of $\vac$ in the expression of $v$ as a linear combination
of the vectors  \eqref{Fock basis}. 

We construct a map 
$$\text{pl}\colon\Gr^0\mbC[z^{-1},z]]\to\cF^{[0]}\, $$
by the following rule.
 Let $H\in\Gr^0\mbC[z^{-1},z]]$. Let 
$$f_1,f_2,f_3,\ldots\in\mbC[z^{-1},z]]$$ be a basis in $H$ of the form 
$f_i(z)=z^{-i}(1+o(1))$. Let 
$$
\text{pl}(H)=f_1\wedge f_2\wedge f_3\wedge\ldots.
$$
The infinite
wedge product is defined by  picking a single monomial summand in each 
$f_i$ in such a way that the summand is $z^{-i}$ for all but finitely many indices $i$ (and summing over possible such choices). The resulting vector
in $\cF^{[0]}$ is easily seen to be independent of the
basis choice $\{f_i\}$.

\subsubsection{Tau-functions from the infinite Grassmanian}\label{subsubsection:tau-functions from Grassmanian}

For any integer $k$, define the operator 
$$\psi_k\colon\cF\to\cF\, , \ \ \ \ \
\psi_k\left(z^{a_1}\wedge z^{a_2}\wedge\ldots\right) = z^k\wedge z^{a_1}\wedge z^{a_2}\wedge\ldots\, .
$$
The operator $\psi_k$ 
increases the charge by $1$. Denote by $\psi_k^*$ the 
associated contraction operator,
$$
\psi^*_k\colon z^{a_1}\wedge z^{a_2}\wedge\ldots\mapsto
\begin{cases}
(-1)^{j-1}z^{a_1}\wedge\ldots\wedge\widehat{z^{a_j}}\wedge\ldots,&\text{if there exists $j$ such that $a_j=k$},\\
0,&\text{otherwise}.
\end{cases}
$$
The hat above denotes an omitted element in the wedge product. 
The operator $\psi_k^*$ decreases the charge by $1$. 

For $n\ge 1$, define $\alpha_n=\sum_{i\in\Z}\psi_i\psi_{i+n}^*$. These operators
$\alpha_n$ do not change the charge and therefore 
leave invariant the space $\cF^{[0]}$. 
The operator $\Gamma$ is defined by 
$$
\Gamma=\exp\left(\sum_{n\ge 1}T_n\alpha_n\right).
$$

For $H\in\Gr^0\mbC[z^{-1},z]]$, define a series $\tau_H(T_1,T_2,T_3,\ldots)$ by
$$
\tau_H(T_1,T_2,T_3,\ldots)=\Big<0\, \Big|\, \Gamma \big(\text{pl}(H)\big)\Big>\, .
$$
The series $\tau_H$ is a tau-function of the KP hierarchy, see \cite{DJKM83}.

\subsubsection{ The partition function
$\exp(F^c)$ as a point in the infinite Grassmanian}\label{subsubsection:point}

For $i\geq 1$,
 define the Laurent series $f_i$:
$$
f_i=
\begin{cases}
z^{-i}{\mathcal{A}}(-z),&\text{if $i$ is odd},\\
-z^{-i}{\mathcal{B}}(-z),&\text{if $i$ is even}.
\end{cases}
$$
Let $H_{\mathcal{A},\mathcal{B}}\subset\mbC[z^{-1},z]]$ be the subspace spanned by the Laurent series $f_i$. The following result is proven by
Kac and Schwarz \cite{KS91}: 
\begin{equation}\label{kqqk}
\left.\exp(F^c)\right|_{t_i=(2i+1)!!T_{2i+1}}=\tau_{H_{\mathcal{A},\mathcal{B}}}\, .
\end{equation}
In fact, \eqref{kqqk} is essentially equivalent to \eqref{eq:determinantal formula},
see \cite[Lemma 4.2]{Kon92}.

\section{The moduli space of Riemann surfaces with boundary}\label{section:open curves}
\subsection{Overview}
The series $\AAA$ also appears in the intersection theory 
of the moduli space of Riemann surfaces with boundary 
(often viewed, with the boundary removed, as open Riemann surfaces).
We recall the basics of the moduli of Riemann surfaces with boundary in Section~\ref{subsection:open moduli}
following \cite{PST}.
In Section~\ref{subsection:open KdV and open Virasoro},
 we review two equivalent conjectural descriptions of 
the descendent theory: the open KdV and the open Virasoro equations. 
An explicit formula for the partition
function of the open theory is discussed in 
 Section~\ref{subsection:explicit formula}. 
The hypergeometric series $\AAA$ plays a basic role in
the formula.

\subsection{Moduli of Riemann surfaces with boundary}\label{subsection:open moduli}


Let $\Delta\in\mbC$ be the open unit disk, and let $\overline\Delta$ be its closure. An extendable embedding of the open disk $\Delta$ in a compact Riemann surface $f\colon\Delta\to C$ is a holomorphic map which
 can be extended to a holomorphic embedding of an open neighborhood of $\overline\Delta$. Two extendable embeddings are disjoint, if the images of $\overline\Delta$ are disjoint.

A Riemann surface with boundary $(X,\d X)$ is obtained by removing a  
finite positive number of disjoint extendable open disks from a connected compact Riemann surface. A compact Riemann surface is {\em not} viewed
here as Riemann surface with boundary.

To a Riemann surface with boundary $(X,\d X)$, we can canonically 
construct the double via the Schwartz reflection through the boundary. 
The double $D(X,\d X)$ of $(X,\d X)$ is a compact Riemann surface. 
The doubled genus of $(X,\d X)$ is defined to be the usual genus of $D(X,\d X)$.

On a Riemann surface with boundary $(X,\d X)$, we consider 
two types of marked points. The markings of interior type 
are points of $X\backslash\d X$. The markings of boundary type 
are points of $\d X$. Let $\mcM_{g,k,l}$ denote the moduli space of 
Riemann surfaces with boundary of doubled genus $g$ with
$k$ distinct boundary markings and $l$ distinct interior markings. 
The moduli space $\mcM_{g,k,l}$ is defined to be empty
 unless the stability condition $$2g-2+k+2l>0$$ is satisfied. 
The moduli space $\mcM_{g,k,l}$ is a real orbifold of 
real dimension $3g-3+k+2l$. 

The cotangent line
classes $\psi_i\in H^2(\mcM_{g,k,l},\mathbb{Q})$ are defined 
(as before)
as the first Chern classes of the cotangent line bundles 
associated to the interior markings. 
In \cite{PST}, cotangent lines at the boundary points are not
considered. Open intersection numbers are defined by
\begin{gather}\label{eq: open intersections}
\left<\tau_{a_1}\tau_{a_2}\ldots\tau_{a_l}\sigma^k\right>^o_g=\int_{\oM_{g,k,l}}\psi_1^{a_1}\psi_2^{a_2}\ldots\psi_l^{a_l}\, .
\end{gather}
To rigorously define the right-hand side of~\eqref{eq: open intersections}, at least three significant steps must be taken:
\begin{itemize}

\item A natural compactification $\mcM_{g,k,l}\subset\oM_{g,k,l}$ must be constructed. Candidates for $\oM_{g,k,l}$ are themselves real orbifolds with boundary $\d\oM_{g,k,l}$;

\item For integration over $\oM_{g,k,l}$ to be well-defined, boundary conditions of the integrand along $\d\oM_{g,k,l}$ must be specified;

\item Orientation issues should be resolved, since the moduli space $\mcM_{g,k,l}$ is in general non-orientable.  

\end{itemize}
All three steps are completed in genus 0
in \cite{PST}. The higher genus constructions
will appear in upcoming work of Solomon and Tessler \cite{ST}.

 We introduce formal variables $t_0,t_1,t_2,\ldots$ and $s$. The generating series $F^o$ is defined by
$$
F^o(t_0,t_1,\ldots,s)=\sum_{\substack{g,k,l\ge 0\\2g-2+k+2l>0}}\frac{1}{k!l!}\sum_{a_1,\ldots,a_l\ge 0}\<\tau_{a_1}\ldots\tau_{a_l}\sigma^k\>^o_gt_{a_1}\ldots t_{a_l} s^k.
$$ 
The series $F^o$ is the {\em open potential}.


\subsection{Open KdV and open Virasoro equations}\label{subsection:open KdV and open Virasoro}

\subsubsection{Constraints}
 KdV and Virasoro type constraints for the open intersection numbers~\eqref{eq: open intersections} were conjectured in \cite{PST}
for all genera (and proven in  genus $0$). 
The following initial condition follows easily from the definitions:
\begin{gather}\label{eq: initial conditions}
\left.F^o\right|_{t_{i\ge 1}=0}=\frac{s^3}{6}+t_0 s\, .
\end{gather}

\subsubsection{Open KdV equations}

The following system of partial differential equations
for a series $$F\in\mathbb{Q}[[t_0,t_1,t_2\ldots,s]]$$
was introduced in \cite{PST}:
\begin{gather}\label{eq:first half}
\frac{2n+1}{2}\frac{\d F}{\d t_n}=\frac{\d F}{\d s}\frac{\d F}{\d t_{n-1}}+\frac{\d^2 F}{\d s\d t_{n-1}}+\frac{1}{2}\frac{\d F}{\d t_0}\frac{\d^2 F^c}{\d t_0\d t_{n-1}}-\frac{1}{4}\frac{\d^3 F^c}{\d t_0^2\d t_{n-1}}\, ,\quad n\ge 1.
\end{gather} 
The above system is called the {\em open KdV equations}.
The open potential $F^o$ was conjectured in \cite{PST}  to be a solution of 
the open
KdV equations. The open KdV equations \eqref{eq:first half}, the initial condition \eqref{eq: initial conditions}, and the potential $F^c$ together uniquely determine the series $F^o$.
However, 
the existence of a such a solution, proven in \cite{Bur14a} is non-trivial.
 We will denote the unique solution by~$\tF^o$.

\subsubsection{Open Virasoro equations}

The classical Virasoro operators $\{L_n\}_{n\ge -1}$ 
which appear in the descendent theory of closed Riemann surfaces are defined as follows:
\begin{multline*}
\hspace{40pt} L_n\ \ =\ \ \sum_{i\ge 0}\frac{(2i+2n+1)!!}{2^{n+1}(2i-1)!!}(t_i-\delta_{i,1})\frac{\d}{\d t_{i+n}}\\
+\frac{1}{2}\sum_{i=0}^{n-1}\frac{(2i+1)!!(2n-2i-1)!!}{2^{n+1}}\frac{\d^2}{\d t_i\d t_{n-1-i}}
\ \ +\delta_{n,-1}\frac{t_0^2}{2}+\delta_{n,0}\frac{1}{16}\, .
\end{multline*}
The following modified operators, 
\begin{gather*}
\cL_n=L_n+\left(s\frac{\d^{n+1}}{\d s^{n+1}}+\frac{3n+3}{4}\frac{\d^n}{\d s^n}\right),\quad n\ge -1,
\end{gather*}
were introduced in \cite{PST} and conjectured to constrain the partition function,
\begin{gather}\label{eq:open Virasoro equations}
\cL_n\exp(F^o+F^c)=0,\quad n\ge -1\, .
\end{gather}
The  equations \eqref{eq:open Virasoro equations} 
are called the {\em open Virasoro equations}. Equations~\eqref{eq:open Virasoro equations}, 
the initial condition 
$$F^o|_{t_{i\geq0}}=\frac{s^3}{6}\, ,$$ 
and the potential $F^c$ together uniquely determine the series $F^o$.

The power series $\tF^o$ is proven in \cite{Bur14a} to satisfy the open Virasoro equations. 
Hence, the two conjectural descriptions of the open descendent theory, given by the 
open KdV equations and the open Virasoro equations, are equivalent.


\subsection{Formula for the open potential}\label{subsection:explicit formula}
Let $G_z$ be the shift operator which acts on
a series $f(t_0,t_1,\ldots)\in\mbC[[t_0,t_1,\ldots]]$ by
\begin{gather*}
G_z(f)\big(t_0,t_1,t_2,\ldots\big)=f\left(t_0-\frac{k_0}{z},t_1-\frac{k_1}{z^3},t_2-\frac{k_2}{z^5},\ldots\right)\, ,
\end{gather*}
where $k_n=(2n-1)!!$ and, by definition,  $(-1)!!=1$.

Define the numbers $\{ d_i\}_{i\geq 0}$ by
\begin{gather*}
d_n=\sum_{i=0}^n 3^i|a_{n-i}|\prod_{k=1}^i\left(n+\half-k\right),
\end{gather*}
where $a_{n-i}$ is a coefficient of $\mathcal{A}$, see Section 1.
The series $D(x)$, defined by 
$$D(x)=1+\sum_{i\ge 1}d_i x^{3i}\, ,$$
has the following equivalent description: $D(x)$
is the unique 
 series solution of the differential equation
$$
\left(-x^4\frac{\d}{\d x}-\frac{3}{2}x^3+1\right)C(x)=\mathcal{A}(-x).
$$ 

The formula of \cite{Bur14b} expressing the open descendent theory
in terms of the closed also requires the series
$$
\xi(t,s;z)=   \frac{s}{2}z^2 +\sum_{i\ge 0}\frac{t_i}{(2i+1)!!}z^{2i+1}\, .
$$
\begin{theorem} [Buryak] 
We have
\begin{gather}\label{eq:explicit formula}
\exp(\tF^o)=\Coef_{z^0}\left[D\left(z^{-1}\right)\frac{G_z(\exp\left(F^c)\right)}{\exp(F^c)}\exp(\xi)\right].
\end{gather}
\end{theorem}
The product
$D\left(z^{-1}\right)\frac{G_z(\exp\left(F^c)\right)}{\exp(F^c)}$ is a 
series in $z^{-1}$. On the other hand, $\exp(\xi)$ is a series in $z$. 
In general, the multiplication of two such series
may not be well-defined. 
In our case, the issue is resolved as follows. 
We introduce a grading in the ring $\mbC[[t_0,t_1,t_2,\ldots,s]]$ 
assigning to $t_i$ the degree $2i+1$ and to $s$ the degree $2$. 
Since the degree of the coefficient of~$z^i$ in $\exp(\xi)$ grows as
$i$ grows, the product in the square brackets is well defined.

In fact, a more general statement is proven in~\cite{Bur14b}.  
A natural way to include variables 
$$s_1,s_2,s_3, \ldots$$ in the series $\tF^o$
was proposed in \cite{Bur14a}.
The new variables $s_i$ may be viewed as descendants of the boundary marked points. The extended power series is denoted by $\tF^{o,ext}$. 
In~\cite{Bur14b}, a formula similar to~\eqref{eq:explicit formula} is proven for the extended series $\tF^{o,ext}$.
The series
$$
\left.\exp(\widetilde F^{o,ext}+F^c)\right|_{\substack{\, t_i=(2i+1)!!T_{2i+1}\ \ \ \\s_i=2^{i+1}(i+1)!T_{2i+2}}}
$$
is proven to be 
a tau-function of the KP hierarchy
in \cite{Ale15} using the formula of \cite{Bur14b}.
\begin{question}
  Is there a system of tautological relations involving $\AAA$ and
  $\BBB$ in the cohomology of the moduli spaces $\oM_{g,k,l}$ parallel
  to Pixton's?
\end{question}

\section{Cohomological field theories and Witten's 3-spin
  class} \label{sec:3spin}

\subsection{Overview}
We describe here how the hypergeometric series $\AAA$ and $\BBB$
appear in the study of Witten's 3-spin class via cohomological field
theories.

In Section~\ref{sec:3spin:CohFTs}, we recall the basics of 
cohomological field theories and the Givental-Teleman classification
\cite{Giv01,Tel} in the  semisimple case. In
Section~\ref{sec:3spin:3spin}, we consider the example of Witten's
$r$-spin class and see how the $\AAA$ and $\BBB$ series appear in the
$R$-matrix of the 3-spin theory. Finally, in
Sections ~\ref{sec:3spin:frob} and \ref{xxzz}, we see how the Airy differential equation
is directly related to the flatness equation for the Dubrovin
connection of the 3-spin theory.

\subsection{Cohomological field theories} \label{sec:3spin:CohFTs}

Cohomological field theories (CohFTs) were first defined by Kontsevich
and Manin \cite{KM97} in order to place the axioms of Gromov-Witten
theory in an algebraic structure.

An $n$-dimensional CohFT, defined on an $n$-dimensional vector space
$V$ together with a nondegenerate bilinear form 
$$\eta: V \times V \rightarrow V$$ and unit vector
$\mathbf 1 \in V$, is a collection of $S_n$-symmetric, multilinear maps
\begin{equation*}
  \Omega_{g, n}: V^n \to H^* (\oM_{g, n},\mathbb{Q}),
\end{equation*}
for every $g$ and $n$ satisfying $2g - 2 + n > 0$, for which
the
following two properties hold:
\begin{itemize}
\item[$\bullet$]{\em Splitting.\  } 
The pull-back of $\Omega_{g, n}(v_1, \dotsc, v_n)$ via a
  glueing map $$\prod_i \oM_{g_i, n_i} \to \oM_{g, n}$$ is the product
  of the $\Omega_{g_i, n_i}$ corresponding to the components, with
  arguments $v_1, \dotsc, v_n$ at the preimages of the marked points
  and the symmetric bivector $\eta^{-1}$ at the points which are glued
  together.
\vspace{6pt}
\item[$\bullet$]{\em Unit.\  }  Let $\pi: \oM_{g, n + 1} \to \oM_{g, n}$ be the
  forgetful map. Then
  \begin{equation*}
    \pi^* \Omega_{g, n}(v_1, \dotsc, v_n)
    = \Omega_{g, n + 1}(v_1, \dotsc, v_n, \mathbf 1).
  \end{equation*}
  Additionally, $\Omega_{0, 3}(v, w, \mathbf 1) = \eta(v, w)$.
\end{itemize}

\smallskip

If all classes of $\Omega$ are of cohomological degree $0$ on the
moduli spaces of curves, $\Omega$ is a \emph{topological field theory}
(TQFT). Then, $\Omega$ is uniquely determined from $\Omega_{0, 3}$ and
$\eta$ by calculating $\Omega_{g, n}$ -- a multiple of the fundamental
class -- at a maximally degenerate curve with $2g - 2 + n$ rational
irreducible components, each with three special points.

The study of CohFTs is motivated by Gromov-Witten theory. To any
nonsingular
projective variety $X$, we can associate
 a CohFT based on the cohomology
ring $H^*(X,\mathbb{Q})$ together with the Poincar\'e pairing
as the bilinear form and with the
fundamental class as the unit vector
(if $X$ has cohomology in odd degree, the
$S_n$-symmetry hypothesis of a CohFT
must be replaced by appropriate skew-symmetry). Let $\oM_{g, n}(X)$ denote the
space of stable maps to $X$, $\pi$ the projection to $\oM_{g, n}$, 
$\ev_i$ the $i$th evaluation map, and $[\oM_{g, n}(X)]^{vir}$ the
virtual fundamental class. We define
\begin{equation} \label{eq:GWCohFT}
  \Omega_{g, n}(v_1, \dotsc, v_n)
  = \pi_* \left(\prod_{i = 1}^n \ev_i^*(v_i) \cap [\oM_{g, n}(X)]^{vir}\right)
\end{equation}
under convergence conditions (required for
the implicit sum over curve classes $\beta\in H_2(X,\mathbb{Z})$
to be well-defined).

The tensors $\Omega_{0, 3}$ and $\eta$ can be used to define a product $\star$,
the \emph{quantum product}, on $V$ via
\begin{equation}
\label{t66}
  \Omega_{0, 3}(a, b, c) = \eta(a \star b, c),
\end{equation}
making $(V, \eta, \star)$ a Frobenius algebra. 
The CohFT $\Omega$
is \emph{semisimple} if the algebra $V$ has no nilpotent
elements, or equivalently, an orthogonal basis of
idempotent elements defined over $\mbC$.

Semisimple CohFTs have been classified by Teleman \cite{Tel} by generalizing a
conjecture of Givental \cite{Giv01}. To reconstruct a CohFT $\Omega$ from the TQFT
$\omega$ defined by the degree $0$ part, a unique
endomorphism valued matrix of power series 
\begin{equation*}
  R = 1 + R_1 z + R_2 z^2 + \dotsb \ \in \End(V)[[z]]  
\end{equation*}
must be specified.  Given $R$, there is a
concrete formula for $\Omega$ in terms of tautological classes similar
to the form of Pixton's relations as in Appendix~\ref{Ssec:relations}.

\subsection{Witten's 3-spin class} \label{sec:3spin:3spin}

For every integer $r \ge 2$, there is a beautiful CohFT obtained from
Witten's $r$-spin class. We review the basic properties of the
construction.

Let $V$ be an $(r - 1)$-dimensional $\Q$-vector space with basis $e_0
, \dotsc, e_{r - 2}$, bilinear form
\begin{equation*}
  \eta(e_a, e_b) = \delta_{a + b, r - 2}\, ,
\end{equation*}
and unit vector $\mathbf 1 = e_0$. Witten's $r$-spin theory provides a
family of classes 
$$W_{g, n} (a_1 , \dotsc, a_n) \in H^*(\oM_{g, n},
\Q)\, $$
 for $a_1 , \dotsc, a_n \in \{0,\dotsc , r - 2\}$. These define a
CohFT on $V$ by setting
\begin{equation*}
  W_{g, n} (e_{a_1} , \dotsc, e_{a_n}) = W_{g, n} (a_1 , \dotsc, a_n)
\end{equation*}
and extending multilinearly.

Witten's class is homogeneous of (complex) degree
\begin{equation*} \label{eq:degwitten}
 \frac{(r-2)(g - 1) + \sum_i a_i}r
\end{equation*}
and vanishes if the degree formula fails to yield an integer.

In genus 0, the construction of Witten's class was first carried
out by Witten \cite{Wit93} using $r$-spin structures ($r^{th}$ roots of
the canonical bundle). In higher genus, there are by now several
constructions. Algebraic approaches have been found by Polishchuk-Vaintrob \cite{PV01} (later 
simplified by Chiodo \cite{Chi06}) and
Chang-Li-Li \cite{CLL}. Analytic constructions  by Mochizuki \cite{Moc06} and 
Fan-Jarvis-Ruan \cite{FJR07} are also available. 
The equivalence of
these constructions is established in \cite{PPZ13}.

\smallskip

The CohFT determined by Witten's $r$-spin class is not semisimple.
For example, for $r = 3$, the
quantum product
$$e_1 \star e_1=0$$ vanishes.  However, via a shift on
the Frobenius manifold, Witten's class can be modified to be a
semisimple CohFT called the \emph{shifted $r$-spin Witten's
  class}. In the case of $r = 3$, the shift depends on one parameter
$\phi$ and, in the new quantum-product,
$$e_1 \star e_1 = \phi e_0\, .$$ The modification destroys
the homogeneity property of Witten's class. The shifted Witten's class
is supported in cohomological degrees at most the degree of Witten's
class.


The $R$-matrix of the CohFT of Witten's $3$-spin class was
calculated explicitly in \cite{PPZ13}. Here, the $\AAA$ and
$\BBB$  hypergeometric series appear. When written
in the basis $\{e_0, e_1\}$, we have
\begin{equation} \label{eq:3spinR}
  R(6 z) =
  \begin{pmatrix}
    -\BBB^{even}(-z \phi^{-3/2}) & -\phi^{1/2} \BBB^{odd}(-z \phi^{-3/2}) \\
    -\phi^{-1/2} \AAA^{odd}(-z \phi^{-3/2}) & \AAA^{even}(-z \phi^{-3/2})
  \end{pmatrix},
\end{equation}
where we have used the superscripts $even$ and $odd$ to denote the
even and odd degree part of a series.

Because the CohFT defined by Witten's class is homogeneous there is a
recursive procedure explicitly by described by Givental and Teleman for
calculating $R$. The matrix \eqref{eq:3spinR} satisfies the recursion
and therefore is the correct $R$-matrix for Witten's 3-spin class.
In the next Sections, we will see a direct connection between the
$R$-matrix and the Airy differential equation.

Applying the Givental-Teleman reconstruction to the $r$-spin CohFT
gives an alternative expression for shifted Witten's class in terms of
tautological classes. The formula of the reconstructed CohFT 
has also terms of cohomological degree higher than
\eqref{eq:degwitten}. The necessary cancellation implies nontrivial
relations between tautological classes. In \cite{PPZ13}, 
these relations for $r=3$ are directly shown to be equivalent to Pixton's
relations. For higher $r$, the relations are studied in \cite{PPZ15}.

\subsection{Frobenius manifolds} \label{sec:3spin:frob}

An $n$-dimensional vector space $V$ may be viewed as a manifold covered
by a chart with coordinates 
$t_1, \ldots, t_n$ corresponding
to a basis $e_1, \ldots, e_n$ of $V$. The global vector fields
$\frac\d{\d t_\mu}$ may be used to identify each tangent space with
$V$.  If $V$ is the vector space of a CohFT, all the tangent spaces of the
manifold $V$ are equipped with the (constant) metric $\eta$. The full genus 0
potential associated to $\Omega$,
$$\Phi(\xi)=\sum_{n\geq 3} \frac{1}{n!}\int_{\overline{\mathcal{M}}_{0,n}} 
\Omega_{0,n}(\xi, \ldots, \xi)\, , \ \ \ \ \ \xi\in V\, ,$$
satisfies the WDVV equations and thus determines 
a deformed quantum product $\star$ on the tangent space at
 $\xi\in V$,
$$\frac{\partial^3\Phi}{\partial t_a \partial t_b \partial t_c}\Big|_\xi
 =
\eta \left(\frac{\partial}{\partial t_a}\star 
\frac{\partial}{\partial t_b}, \frac{\partial}{\partial t_c}\right)\, .$$

Together, the above constructions  endow $V$ with the structure of a  
\emph{Frobenius manifold}. At the origin of $V$ with
coordinates $$t_1 = \ldots = t_n = 0\, ,$$
the quantum product $\star$ agrees with the earlier definition
\eqref{t66}. Frobenius
manifolds have been introduced by Dubrovin and the full definitions can
be found in his monograph \cite{Dub96}.

The coordinates $t_\mu$ are called \emph{flat coordinates}. If the
multiplication $\star$ on the tangent
space of $\xi\in V$ is semisimple,
an alternative
set of \emph{canonical coordinates} $u_i$ are defined
in a neighborhood of $\xi$. These are defined up to
additive constants and reordering by requiring that the corresponding
vector fields $\frac\d{\d u_i}$ form a basis of orthogonal idempotents
at each point where they are defined. 
Let $\mathbf u$ be the diagonal
matrix with the functions $u_i$
along the diagonal,
$$\mathbf{u}=  \left( \begin{array}{ccc}
 u_1&           & \\
    &   \ddots &           \\
        &  &       u_n \end{array}\right)\, .$$
 Let $\Psi$ be the base
change matrix from the basis of vector fields $\frac\d{\d t_\mu}$ to
the basis of \emph{normalized idempotents} given by
$$\eta\left(\frac\d{\d u_i}, \frac\d{\d u_i}\right)^{-1/2} \frac\d{\d u_i}\,  .$$

Up to constants of integration, the {\em $R$-matrix} written in the basis of
normalized idempotents is uniquely determined by the property that the
product
\begin{equation} \label{eq:RS}
  S = \Psi R e^{\mathbf u/z}
\end{equation}
is a matrix of asymptotic fundamental solutions to the
\emph{flatness equation}
\begin{equation*}
  z \frac\d{\d t_\mu}S = \frac\d{\d t_\mu} \star S\, .
\end{equation*}
The name of the equation stems from the fact that it characterizes
flat vector fields for the \emph{Dubrovin (projective) connection}
$\nabla_z$ defined by
\begin{equation*}
  \nabla_{z, X} = z\nabla_X - X \star
\end{equation*}
for any vector field $X$. Here, $\nabla$ is the Levi-Civita connection
corresponding to the metric $\eta$.

\subsection{Witten's 3-spin theory}
\label{xxzz}

The 3-spin CohFT determines a 2-dimensional Frobenius manifold with
potential
\begin{equation*}
  \frac 12 t_0^2 t_1 + \frac 1{72} t_1^4
\end{equation*}
and metric
\begin{equation*}
  \eta =
  \begin{pmatrix}
    0 & 1 \\
    1 & 0
  \end{pmatrix}\,.
\end{equation*}
Therefore, $\frac\d{\d t_0}$ is the unit for the quantum product, and 
\begin{equation*}
  \frac\d{\d t_1} \star \frac\d{\d t_1} = \phi \frac\d{\d t_0}\,,
\end{equation*}
where $\phi = \frac {t_1}3$.

Let $F(x)$ be the versal deformation of the $A_2$-singularity $x^3
= 0$,
\begin{equation*}
  F(x) = x^3 - t_1 x + t_0\, .
\end{equation*}
We can identify the quantum product with the multiplication in the
Milnor ring
\begin{equation*}
  \mbC[t_0, t_1][x] / F'(x)
\end{equation*}
under the identifications $\frac\d{\d t_0} \mapsto 1$ and $\frac\d{\d
  t_1} \mapsto -x$.

The idempotents and canonical coordinates correspond to the critical
points $\pm \sqrt{\phi}$ of $F(x)$. Canonical coordinates are given by
the associated critical values. Explicitly, the  change of basis 
$\Psi$ and the matrix $\mathbf u$ of canonical coordinates are given
by
\begin{align*}
  \Psi =
  \begin{pmatrix}
    \frac {-\sqrt{\phi}}{\sqrt{\Delta_+}} & \frac {\sqrt{\phi}}{\sqrt{\Delta_-}} \\
    \frac 1{\sqrt{\Delta_+}} & \frac 1{\sqrt{\Delta_-}}
  \end{pmatrix}
  , \quad \mathbf u =
  \begin{pmatrix}
    -2 \phi^{3/2} + t_0 & 0 \\
    0 & 2 \phi^{3/2} + t_0
  \end{pmatrix},
\end{align*}
where $\Delta_{\pm} = \mp 2\sqrt{\phi}$, and choices of roots
$\sqrt{\Delta_\pm}$ have been made.

The flatness equations may be written explicitly as the system
\begin{equation} \label{eq:qDE3spin}
  \begin{aligned}
    z \frac \d{\d t_0} S_\pm^0 =& S_\pm^0\, ,\ \ \ &  z \frac \d{\d t_1} S_\pm^1 =& S_\pm^0\, , \\
    z \frac \d{\d t_0} S_\pm^1 =& S_\pm^1\, ,\ \ \ & z \frac \d{\d t_1}
    S_\pm^0 =& \phi S_\pm^1\,,
  \end{aligned}
\end{equation}
where the upper indices $0$ and $1$ stand for the vector components in
the basis $\{e_0, e_1\}$ and the lower index distinguishes two linear
independent solutions. Combining these equations we see that $S_\pm^1$
satisfies the Airy differential equation
\begin{equation} \label{tqq2}
  \left(z\frac \d{\d t_1}\right)^2 S_\pm^1 = \phi S_\pm^1\, .
\end{equation}
Up to recaling $t_1$, the differential equation \eqref{tqq2} is equivalent to 
\eqref{tqq1}.

The solutions $S_\pm^\mu$ are given by the asymptotic expansion for $z
\to 0$ of the complex contour integrals
\begin{equation}\label{rred2}
  \left(\frac{2\pi z}{3}\right)^{-\frac{1}{2}} \int_{\Gamma_\pm} e^{F(x)/z} (-x)^{1 - \mu} \mathrm dx\,
\end{equation}
defined for almost all $(t, z)$. 
Here, the Lefschetz thimbles
$\Gamma_\pm$ correspond to the two critical points 
$$p_\pm =
\pm\sqrt{\phi}$$ of $F(x)/z$ and are chosen as follows (see
also \cite{WW}). 
Suppose
the critical values of $F(x)/z$ have different imaginary values. We
consider $\Re(F/z)$ as function in the real and imaginary part of
$x$. The cycle $\Gamma_\pm$ is the union of two integral curves of the
vector field $-\nabla\Re(F/z)$ arriving at $p_\pm$ at time
$-\infty$. By construction, the real part of $F/z$ decreases fast
enough when moving along $\Gamma_\pm$ from $p_\pm$ such that the
contour integrals converge absolutely. By differentiating
under the integral, the
contour integrals for any choice of cycle are easily seen to give solutions to
\eqref{eq:qDE3spin}.

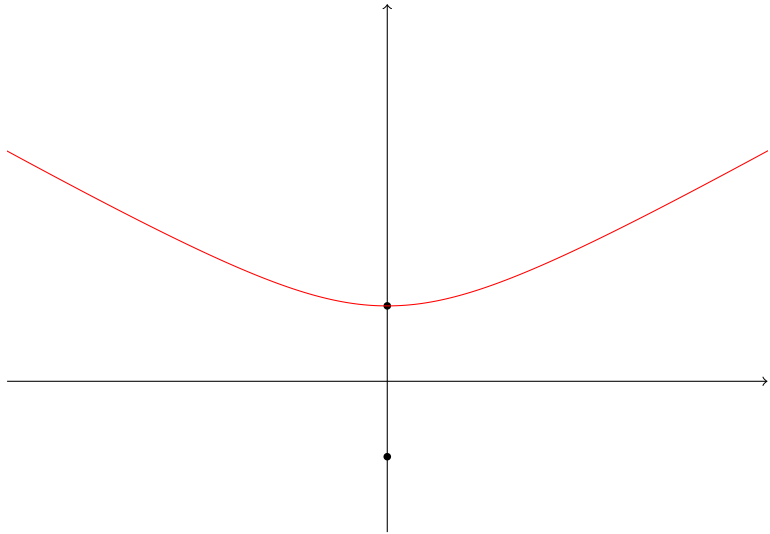
\begin{figure}
\begin{tikzpicture}
\clip (-5, -2) rectangle (5, 5);
\draw[->] (-5, 0) -- (5, 0);
\draw[->] (0, -2) -- (0, 5);
\fill (0, 1) circle (.5mm);
\fill (0, -1) circle (.5mm);
\draw[red] (0, 1.000000) -- (0.050000, 1.000000) -- (0.099984, 1.001250) -- (0.149949, 1.003120)
  -- (0.199877, 1.005815) -- (0.249753, 1.009324) -- (0.299567, 1.013636) -- (0.349306, 1.018737)
  -- (0.398960, 1.024611) -- (0.448519, 1.031239) -- (0.497974, 1.038601) -- (0.547317, 1.046677)
  -- (0.596543, 1.055442) -- (0.645645, 1.064875) -- (0.694619, 1.074951) -- (0.743462, 1.085645)
  -- (0.792171, 1.096934) -- (0.840745, 1.108792) -- (0.889182, 1.121195) -- (0.937482, 1.134120)
  -- (0.985647, 1.147544) -- (1.033676, 1.161443) -- (1.081572, 1.175796) -- (1.129336, 1.190581)
  -- (1.176970, 1.205778) -- (1.224478, 1.221367) -- (1.271862, 1.237329) -- (1.319124, 1.253647)
  -- (1.366269, 1.270302) -- (1.413298, 1.287279) -- (1.460216, 1.304562) -- (1.507026, 1.322135)
  -- (1.553731, 1.339986) -- (1.600334, 1.358101) -- (1.646839, 1.376466) -- (1.693249, 1.395071)
  -- (1.739567, 1.413903) -- (1.785796, 1.432952) -- (1.831940, 1.452208) -- (1.878000, 1.471661)
  -- (1.923981, 1.491301) -- (1.969885, 1.511122) -- (2.015715, 1.531113) -- (2.061473, 1.551268)
  -- (2.107161, 1.571578) -- (2.152784, 1.592039) -- (2.198341, 1.612641) -- (2.243837, 1.633381)
  -- (2.289274, 1.654250) -- (2.334652, 1.675245) -- (2.379975, 1.696360) -- (2.425244, 1.717590)
  -- (2.470462, 1.738929) -- (2.515630, 1.760374) -- (2.560749, 1.781920) -- (2.605822, 1.803563)
  -- (2.650850, 1.825300) -- (2.695835, 1.847126) -- (2.740778, 1.869038) -- (2.785680, 1.891033)
  -- (2.830543, 1.913108) -- (2.875368, 1.935259) -- (2.920157, 1.957484) -- (2.964911, 1.979781)
  -- (3.009630, 2.002145) -- (3.054317, 2.024576) -- (3.098971, 2.047070) -- (3.143594, 2.069625)
  -- (3.188188, 2.092240) -- (3.232752, 2.114912) -- (3.277289, 2.137639) -- (3.321798, 2.160419)
  -- (3.366281, 2.183251) -- (3.410738, 2.206132) -- (3.455170, 2.229062) -- (3.499579, 2.252037)
  -- (3.543964, 2.275058) -- (3.588326, 2.298123) -- (3.632667, 2.321230) -- (3.676986, 2.344377)
  -- (3.721284, 2.367564) -- (3.765563, 2.390789) -- (3.809822, 2.414052) -- (3.854062, 2.437350)
  -- (3.898284, 2.460684) -- (3.942487, 2.484051) -- (3.986674, 2.507451) -- (4.030843, 2.530883)
  -- (4.074996, 2.554346) -- (4.119133, 2.577840) -- (4.163254, 2.601362) -- (4.207361, 2.624913)
  -- (4.251452, 2.648491) -- (4.295529, 2.672096) -- (4.339592, 2.695727) -- (4.383642, 2.719384)
  -- (4.427678, 2.743065) -- (4.471702, 2.766771) -- (4.515712, 2.790499) -- (4.559711, 2.814251)
  -- (4.603697, 2.838025) -- (4.647672, 2.861820) -- (4.691636, 2.885636) -- (4.735588, 2.909473)
  -- (4.779530, 2.933330) -- (4.823461, 2.957206) -- (4.867381, 2.981101) -- (4.911292, 3.005015)
  -- (4.955192, 3.028946) -- (4.999084, 3.052896) -- (5.042965, 3.076862) -- (5.086838, 3.100845); 
\draw[red] (0, 1.000000) -- (-0.050000, 1.000000) -- (-0.099984, 1.001250) -- (-0.149949, 1.003120)
  -- (-0.199877, 1.005815) -- (-0.249753, 1.009324) -- (-0.299567, 1.013636) -- (-0.349306, 1.018737)
  -- (-0.398960, 1.024611) -- (-0.448519, 1.031239) -- (-0.497974, 1.038601) -- (-0.547317, 1.046677)
  -- (-0.596543, 1.055442) -- (-0.645645, 1.064875) -- (-0.694619, 1.074951) -- (-0.743462, 1.085645)
  -- (-0.792171, 1.096934) -- (-0.840745, 1.108792) -- (-0.889182, 1.121195) -- (-0.937482, 1.134120)
  -- (-0.985647, 1.147544) -- (-1.033676, 1.161443) -- (-1.081572, 1.175796) -- (-1.129336, 1.190581)
  -- (-1.176970, 1.205778) -- (-1.224478, 1.221367) -- (-1.271862, 1.237329) -- (-1.319124, 1.253647)
  -- (-1.366269, 1.270302) -- (-1.413298, 1.287279) -- (-1.460216, 1.304562) -- (-1.507026, 1.322135)
  -- (-1.553731, 1.339986) -- (-1.600334, 1.358101) -- (-1.646839, 1.376466) -- (-1.693249, 1.395071)
  -- (-1.739567, 1.413903) -- (-1.785796, 1.432952) -- (-1.831940, 1.452208) -- (-1.878000, 1.471661)
  -- (-1.923981, 1.491301) -- (-1.969885, 1.511122) -- (-2.015715, 1.531113) -- (-2.061473, 1.551268)
  -- (-2.107161, 1.571578) -- (-2.152784, 1.592039) -- (-2.198341, 1.612641) -- (-2.243837, 1.633381)
  -- (-2.289274, 1.654250) -- (-2.334652, 1.675245) -- (-2.379975, 1.696360) -- (-2.425244, 1.717590)
  -- (-2.470462, 1.738929) -- (-2.515630, 1.760374) -- (-2.560749, 1.781920) -- (-2.605822, 1.803563)
  -- (-2.650850, 1.825300) -- (-2.695835, 1.847126) -- (-2.740778, 1.869038) -- (-2.785680, 1.891033)
  -- (-2.830543, 1.913108) -- (-2.875368, 1.935259) -- (-2.920157, 1.957484) -- (-2.964911, 1.979781)
  -- (-3.009630, 2.002145) -- (-3.054317, 2.024576) -- (-3.098971, 2.047070) -- (-3.143594, 2.069625)
  -- (-3.188188, 2.092240) -- (-3.232752, 2.114912) -- (-3.277289, 2.137639) -- (-3.321798, 2.160419)
  -- (-3.366281, 2.183251) -- (-3.410738, 2.206132) -- (-3.455170, 2.229062) -- (-3.499579, 2.252037)
  -- (-3.543964, 2.275058) -- (-3.588326, 2.298123) -- (-3.632667, 2.321230) -- (-3.676986, 2.344377)
  -- (-3.721284, 2.367564) -- (-3.765563, 2.390789) -- (-3.809822, 2.414052) -- (-3.854062, 2.437350)
  -- (-3.898284, 2.460684) -- (-3.942487, 2.484051) -- (-3.986674, 2.507451) -- (-4.030843, 2.530883)
  -- (-4.074996, 2.554346) -- (-4.119133, 2.577840) -- (-4.163254, 2.601362) -- (-4.207361, 2.624913)
  -- (-4.251452, 2.648491) -- (-4.295529, 2.672096) -- (-4.339592, 2.695727) -- (-4.383642, 2.719384)
  -- (-4.427678, 2.743065) -- (-4.471702, 2.766771) -- (-4.515712, 2.790499) -- (-4.559711, 2.814251)
  -- (-4.603697, 2.838025) -- (-4.647672, 2.861820) -- (-4.691636, 2.885636) -- (-4.735588, 2.909473)
  -- (-4.779530, 2.933330) -- (-4.823461, 2.957206) -- (-4.867381, 2.981101) -- (-4.911292, 3.005015)
  -- (-4.955192, 3.028946) -- (-4.999084, 3.052896) -- (-5.042965, 3.076862) -- (-5.086838, 3.100845); 
\end{tikzpicture}
\caption{Lefschetz thimble for $i(z^3 + 3z)$ through critical point $i$.}
\end{figure}

The Airy function arises (up to a multiplicative constant) in the
case  $$\mu = 1\, ,\  z = -3i\, ,\  t_0 = 0\, , \ \text{and } t_1 > 0\, $$
for the critical point $\sqrt{\phi}$. Under these conditions, the
Lefschetz thimble can be deformed into the cycle in the integral used
in the defintion \eqref{eq:defairy} of the Airy function.

As in Section~\ref{sec:asymp}, we make the integral look like a
Gaussian integral by translating and scaling the integrand,
\begin{equation*}
\left(\frac{2\pi z}{3}\right)^{-\frac{1}{2}}\int_{\Gamma_\pm} e^{F(x)/z} (-x)^{1 - \mu} \mathrm dx = \frac{e^{u_\pm/z}}{\sqrt{2\pi} \sqrt{\Delta_\pm}} \int_{\widetilde \Gamma_\pm} e^{-\frac{x^2}2 - x^3 \sqrt{-z} (-3\Delta_\pm)^{-3/2}} \left(\frac{-x\sqrt{-z}}{\sqrt{-3\Delta_\pm}} \mp \sqrt{\phi}\right)^{1 - \mu} \mathrm dx\,,
\end{equation*}
where $\widetilde \Gamma_\pm$ is the Lefschetz thimble defined by the new
exponent. Because of boundeness of the integral as $z \to 0$ and the
fact that in this limit the cycle $\tilde \Gamma_\pm$ approaches the
real line, we obtain as in Section~\ref{sec:asymp} an asymptotic
expansion by formally expanding the integrand in $z$ and calculating
the integrals of the individual summands. For $\mu = 1$, we obtain
\begin{equation*}
  S^1_\pm e^{-u_\pm/z} = \frac 1{\sqrt{\Delta_\pm}} \sum_{j = 0}^\infty \frac{(6j)!}{(3j)!(2j)!} \left(\frac z{216\Delta_\pm^3}\right)^j = \frac 1{\sqrt{\Delta_\pm}} \AAA\left(\frac{\mp z}{6\phi^{3/2}}\right)\, ,
\end{equation*}
and, for $\mu = 0$, 
\begin{equation*}
  S^0_\pm e^{-u_\pm/z} = \frac{\mp\sqrt\phi}{\sqrt{\Delta_\pm}} \sum_{j = 0}^\infty \frac{(6j)!}{(3j)!(2j)!} \frac{-6j - 1}{6j - 1} \left(\frac z{216\Delta_\pm^3}\right)^j = \frac{\pm\sqrt\phi}{\sqrt{\Delta_\pm}} \BBB\left(\frac{\mp z}{6\phi^{3/2}}\right)\,.
\end{equation*}
By \eqref{eq:RS}, we obtain the $R$-matrix written in flat coordinates
as the product
\begin{equation*}
  \Psi^{-1} S e^{-\mathbf u/z}\, .
\end{equation*}
By the formula
\begin{equation*}
  \Psi^{-1} =
  \begin{pmatrix}
    \frac 1{\sqrt{\Delta_+}} & \frac{-\sqrt\phi}{\sqrt{\Delta_+}} \\
    \frac 1{\sqrt{\Delta_-}} & \frac{\sqrt\phi}{\sqrt{\Delta_-}}
  \end{pmatrix}\, ,
\end{equation*}
we immediately arrive at \eqref{eq:3spinR}.

\section{Stable maps and quotients with target  $\CP1$}
\subsection{Antecedents}
The $\AAA$ and $\BBB$ series can also be encountered in $R$-matrices of
other CohFTs. Here, we consider the
equivariant Gromov-Witten theory of $\CP1$.

The analysis is motivated by the earlier study of the hypergeometric
function
\begin{equation*}
  \Phi(z, q)
  = \sum_{d = 0}^\infty \prod_{i = 1}^d \frac 1{\lambda - iz} \ \frac{(-1)^d}{d!} \frac{q^d}{z^d}
  = \sum_{d = 0}^\infty q^d \prod_{i = 1}^d \frac 1{(iz - \lambda)iz}
\end{equation*}
which arises naturally in the geometry of the moduli of curves.
Let
$$\mathbb B_d\rightarrow {\mathcal{M}}_{g,d}$$ 
denote the bundle with
fiber $H^0(C, \mathcal O_C(\sum_j p_j)|_{\sum_j p_j})$ over the moduli point
$[C,p_1, \dotsc, p_d]$.
The series $\Phi$ at
$\lambda = 1$ has been used to calculate
\begin{equation*}
  \pi_* c_z^{-1}(\mathbb B_d),
\end{equation*}
where  $c_z^{-1}$ the inverse of
the Chern polynomial in variable $z$ and
$$\pi: \mathcal{C}_{g, n}^d \to \mathcal{M}_{g, n}$$ is the projection from the
$d$-fold universal curve to $\mathcal{M}_{g, n}$.
The study of $\Phi$ was first taken up by Ionel
\cite{Ion05} to prove the nonvanishing of certain coefficients
in tautological relations related to Faber's generation conjecture --
the hypergeometric series $\AAA$ already arises in \cite{Ion05}.
In \cite{PP13}, a further study of $\Phi$ was required as
a part of the proof of the Faber-Zagier relations (and the
series $\BBB$ also emerged). The series $\Phi$ plays a
basic role in the proof \cite{Jan13} of Pixton's
relations in the Chow ring.

The relation to the Gromov-Witten theory of $\CP1$ is as follows. In
\cite{Jan13,PP13}, the bundle $\mathbb B_d$ appeared as a
vertex contribution in a localization calculation for stable quotients
\cite{MOP11} to $\CP1$. There is a proper map from the moduli space of
of stable maps to the moduli space of stable quotients to $\CP1$,
compatible with localization in the sense that the vertex
contributions for stable quotients are sums of localization
contributions for stable maps. Hence, the localization calculation yields the
same results in both spaces. Givental's conjecture on the
reconstruction of CohFTs was motivated by equivariant localization:
the case of the equivariant Gromov-Witten
theory of toric varieties was proven in \cite{Giv01} using 
virtual localization
\cite{GrP}.

\subsection{Frobenius manifold}

The equivariant Gromov-Witten theory $\CP1$ defines a 2 dimensional
CohFT by the equivariant analog of \eqref{eq:GWCohFT}. The genus zero
potential of the corresponding Frobenius manifold records the genus
zero, equivariant primary Gromov-Witten invariants of $\CP1$. It is
given by
\begin{equation*}
  \frac 12 t_0^2 t_1 + \frac\lambda 2 t_0 t_1^2 - \frac{\lambda^2}6 t_1^3 + e^{t_1}
\end{equation*}
and the metric is
\begin{equation*}
  \eta =
  \begin{pmatrix}
    0 & 1 \\
    1 & \lambda
  \end{pmatrix}\,.
\end{equation*}
Here, $\lambda$ is the equivariant
parameter.
The quantum product of the Frobenius manifold coincides with
the equivariant quantum product of $\CP1$, determined by the equation
\begin{equation*}
  H(H - \lambda) = e^{t_1}\, ,
\end{equation*}
after we identify $\frac \d{\d t_0}$ with the fundamental class and
$\frac \d{\d t_1}$ with the hyperplane class $H$.

Let $F$ be the mirror curve, 
\begin{equation*}
  F(x) = e^x + e^{t_1 - x} + \lambda x + t_0\,.
\end{equation*}
Identifying $H$ with $e^{t_1 - x}$, we can interpret the quantum
cohomology ring also as
\begin{equation*}
  \mbC[t_0, e^{t_1}][e^x, e^{-x}] / F'(x)\, .
\end{equation*}
The critical points are
\begin{equation*}
  e^x  
= -\frac\lambda 2 \pm \sqrt{e^{t_1} + \frac{\lambda^2}4}\, .
\end{equation*}
 Setting $\phi = e^{t_1} +\frac{\lambda^2}{4}$, 
we determine the canonical coordinates
\begin{equation*}
  u_\pm = \pm 2\sqrt\phi + \lambda \ln\left(-\frac\lambda 2 \pm \sqrt\phi\right) + t_0
\end{equation*}
and the basis change matrix
\begin{equation*}
  \Psi =
  \begin{pmatrix}
    \frac {-\frac\lambda 2 + \sqrt{\phi}}{\sqrt{\Delta_+}} & \frac {-\frac\lambda 2 - \sqrt{\phi}}{\sqrt{\Delta_-}} \\
    \frac 1{\sqrt{\Delta_+}} & \frac 1{\sqrt{\Delta_-}}
  \end{pmatrix}
\end{equation*}
for $\Delta_\pm = \pm 2\sqrt\phi$.

The flatness equations can be written as
\begin{eqnarray*}
  z \frac \d{\d t_0} S_\pm^0 = S_\pm^0\, ,\ \ \  & & z \frac \d{\d t_1} S_\pm^1 = S_\pm^0 + \lambda S_\pm^1\, , \\
  z \frac \d{\d t_0} S_\pm^1 = S_\pm^1\, ,\ \ \  & & z \frac \d{\d t_1} S_\pm^0 = e^{t_1} S_\pm^1\, ,
\end{eqnarray*}
in the basis $\{1, H\}$. They imply the second order differential
equation
\begin{equation*}
  \left(z\frac \d{\d t_1}\right)^2 S_\pm^1 - \lambda z\frac \d{\d t_1} S_\pm^1
  = e^{t_1} S_\pm^1
\end{equation*}
for $S_\pm^1$. After the variable change
\begin{equation*}
  q = e^{t_1}\, ,  
\end{equation*}
the hypergeometric function $\Phi$ satisfies the same differential
equation.

\subsection{Asymptotic analysis}

Asymptotic solutions to the flatness equation can again be constructed
using asymptotic expansion of contour integrals
\begin{equation*}
  \frac 1{\sqrt{-2\pi z}} \int_{\Gamma_\pm} e^{F(x, t)/z} (e^{t_1 - x} - \lambda)^{1 - \mu} \mathrm dx
\end{equation*}
along Lefschetz thimbles $\Gamma_\pm$ though critical points of $F$.

The series $\Phi$ can be identified with the oscillating integral
solution for $\mu = 0$ and the cycle $\Gamma_-$, up to a factor
independent of $q = e^{t_1}$. To prove this claim, as both
functions satisfy the same second order differential equation in $q$,
we need only check that the functions and their first $q$-derivatives
agree at $q = 0$ up to the same factor.

We first study the limit $q \to 0$ of the integral
\begin{eqnarray*}
  \lim_{q \to 0} \frac 1{\sqrt{-2\pi z}} \int_{\Gamma_-} e^{F(x, t)/z} \mathrm dx
 & = & \frac 1{\sqrt{-2\pi z}} \int_{\Gamma_-} e^{(e^x + \lambda x + t_0)/z} \mathrm dx \\
 & = & \frac{(-z)^{\lambda/z}}{\sqrt{-2\pi z}} \int_{\widehat \Gamma_-} e^{-x + (\lambda \ln(x) + t_0)/z} \frac{\mathrm dx}x\,,
\end{eqnarray*}
where the substitution $x \mapsto \ln(-xz)$ was applied in the second
step. The critical point that the new Lefschetz thimble $\widehat\Gamma_-$
moves through is at $\frac z\lambda$ and, if we assume that this ratio
is positive real, $\widehat\Gamma_-$ coincides with the positive real
axis. So we can rewrite the limit as
\begin{equation*}
  \lim_{q \to 0} \frac 1{\sqrt{-2\pi z}} \int_{\Gamma_-} e^{F(x, t)/z} \mathrm dx
  = \frac{(-z)^{\lambda/z} e^{t_0/z}}{\sqrt{-2\pi z}} \Gamma\left(\frac{-\lambda}{-z}\right),
\end{equation*}
where $\Gamma$ is the Gamma function. Using Stirling's formula
\begin{equation*}
  \ln\Gamma(x) \, \asymp\, x\ln(x) - x - \frac 12 \ln\left(\frac x{2\pi}\right) + \sum_{i = 1}^\infty \frac{B_{2i}}{2i(2i-1)} x^{-2i}\,,
\end{equation*}
where the $B_{2i}$ are the Bernoulli numbers defined by
\begin{equation*}
  \sum_{k = 0}^\infty B_k \frac{x^k}{k!} = \frac x{e^x - 1}\,,
\end{equation*}
we find 
\begin{equation} \label{eq:bermod}
  \lim_{q \to 0} \frac 1{\sqrt{-2\pi z}} \int_{\Gamma_-} e^{F(x, t)/z} \mathrm dx
 \, \asymp\, \frac{e^{u_-|_{q = 0}/z}}{\sqrt{\Delta_-|_{q = 0}}} \exp\left(\sum_{i = 1}^\infty \frac{B_{2i}}{2i(2i - 1)} (z/\lambda)^{2i - 1}\right)\,,
\end{equation}
where $u_-|_{q = 0} = -\lambda + \lambda \ln(-\lambda) + t_0$ and
$\Delta_-|_{q = 0} = -\lambda$. Because of $\Phi(z, 0) = 1$, the
factor we need to multiply $\Phi$ with so that it can coincide with
the oscillating integral is \eqref{eq:bermod}. Geometrically, the
Bernoulli modification is caused by the contribution of the Hodge
bundle in the vertex terms of the localization.

We also need to check that at $q = 0$ the first $q$-derivative of $\Phi$
\begin{equation*}
  \frac{\d \Phi}{\d q} (z, 0) = \frac 1{(z - \lambda)z}
\end{equation*}
coincides with the corresponding oscillating integral up to the factor
\eqref{eq:bermod}. Similarly to above, we compute
\begin{multline*}
  \lim_{q \to 0} \frac{\d}{\d q} \frac 1{\sqrt{-2\pi z}} \int_{\Gamma_-} e^{F(x, t)/z} \mathrm dx
  = \lim_{q \to 0} \frac 1{\sqrt{-2\pi z}} \int_{\Gamma_-} \frac{e^{-x}}z e^{F(x, t)/z} \mathrm dx
  = \frac{(-z)^{\lambda/z - 1} e^{t_0/z}}{z\sqrt{-2\pi z}} \Gamma\left(\frac{-\lambda}{-z} - 1\right) \\
  = \frac 1{z(z - \lambda)} \frac{(-z)^{\lambda/z} e^{t_0/z}}{\sqrt{-2\pi z}} \Gamma\left(\frac{-\lambda}{-z}\right) = \frac 1{z(z - \lambda)} \lim_{q \to 0} \frac 1{\sqrt{-2\pi z}} \int_{\Gamma_-} e^{F(x, t)/z} \mathrm dx\,,
\end{multline*}
completing the proof that $\Phi$ coincides with the oscillating
integral up to the factor \eqref{eq:bermod}.

\medskip

Let us extract the $A$-function from a different limit of the
$R$-matrix. To calculate the asymptotic expansion, we can proceed as in
Section~\ref{xxzz} by formally expanding $\exp(F/z)$ near a critical
point. Let us first expand $F/z$ as
\begin{equation*}
  \frac 1z \left(u_\pm \pm 2\phi \frac{x^2}2 - \lambda \frac{x^3}6 \pm 2\phi \frac{x^4}{4!} - \lambda \frac{x^5}{5!} \pm 2\phi \frac{x^6}{6!} - \dotsb\right)\,,
\end{equation*}
which becomes
\begin{equation*}
  \frac{u_\pm}z -\frac{x^2}2 + \lambda \sqrt{-z} \Delta_\pm^{-3/2} \frac{x^3}6 - (-z) \Delta_\pm^{-1} \frac{x^4}{4!} + \lambda (-z)^{3/2} \Delta_\pm^{-5/2} \frac{x^5}{5!} - (-z)^2 \Delta_\pm^{-2} \frac{x^6}{6!} + \dots
\end{equation*}
after replacing $x \mapsto x \sqrt{-z} \Delta_\pm^{-1/2}$. Notice 
the cubic term in $x$  essentially carries a power of $z/\phi^{3/2}$
whereas the higher order terms  carry powers of
$z/\phi^c$ for $c < \frac 32$. Therefore, if we are only interested in
the coefficients of the $R$-matrix with lowest possible power of
$\phi$, we can ignore all the higher order terms of $F/z$. In the $\mu
= 1$ integral without the higher order terms, we discover the
$\AAA$-series,
\begin{equation*}
  \frac 1{\sqrt{2\pi} \sqrt{\Delta_\pm}} \int_{\tilde\Gamma_\pm} \exp\left(\frac{u_\pm}z -\frac{x^2}2 + \lambda \sqrt{-z} \Delta_\pm^{-3/2} \frac{x^3}6\right) \mathrm dx
 \, \asymp\, \frac{e^{u_\pm}}{\sqrt{\Delta_\pm}} \AAA\left(\frac{\mp z \lambda^2}{8\phi^{3/2}}\right)\,.
\end{equation*}
For $\mu = 0$, we similarly obtain a linear combination of the $\AAA$-
and $\BBB$-series.

Extracting the lowest order terms can also be achieved by taking the following peculiar
limit. After replacing $z$ by $z/\lambda^2$, take the limit
$\lambda \to \infty$ while keeping $\phi$ fixed. This limit does not
make sense immediately for all data: the canonical coordinates $u_\pm$
need to be changed by additive constants and  a flat basis $\{1, H -
\lambda/2\}$ should be used instead of $\{1, H\}$. In the limit, the
$R$-matrix essentially agrees with the $3$-spin $R$-matrix. In
\cite{Jan14}, it is more generally shown how the CohFT of Witten's
$r$-spin class arises as a limit of the equivariant Gromov-Witten
theory of projective $(r - 2)$-space.

\subsection{Universality}
For any Frobenius manifold obtained from a CohFT which
is generically semisimple on the vector space $V$, 
relations in the tautological rings of the moduli spaces
of curves can be found by studying the pole cancellations
required as the limit of the Givental-Teleman formula
to a non-semisimple point is taken. The first comparison
result is the following

\begin{theorem} [Janda] \label{theorem:main}
The tautological relations in $H^*(\oM_{g,n},\mathbb{Q})$ 
obtained from 3-spin and
equivariant $\CP1$ are equivalent.
\end{theorem}

Since the 3-spin relations were proven to be equivalent
to Pixton's relations in \cite{PPZ13}, the
equivariant $\CP1$ relations are also equivalent to Pixton's.
A natural question, perhaps
more approachable than Pixton's completeness conjecture,
is the following.

\begin{question}
  Do all tautological relations obtained from generically semisimple
  CohFTs by taking limits to non-semisimple points lie in Pixton's
  set?
\end{question}

An affirmative answer to Question 3 would imply that
the hypergeometric series $\AAA$ and $\BBB$ are
lurking in the structure of all generically semisimple
CohFTs.

\appendix
\section{Pixton's relations} 
\label{Ssec:relations}
\subsection{Strata algebra}
Let $\cS_{g,n}^*$ be the $\mathbb{Q}$-algebra of
$\kappa$ and $\psi$ classes 
supported on the strata $\overline{\mathcal{M}}_{g,n}$.
A $\mathbb{Q}$-basis of $\cS_{g,n}^*$ is given by isomorphism classes
of pairs $[\Gamma, \gamma]$ where $\Gamma$ is 
stable graph corresponding to a stratum of the
moduli space of curves, 
$$\overline{\mathcal{M}}_\Gamma \rightarrow \overline{\mathcal{M}}_{g,n}\, ,$$
and $\gamma$ is product of $\kappa$ and $\psi$
classes on $\overline{\mathcal{M}}_\Gamma$.
The strata algebra $\cS_{g,n}^*$ is graded by
codimension
$$\cS_{g,n}^* = \bigoplus_{d=0}^{3g-3+n} \cS^d_{g,n}\ $$
and carries a product for which
the natural push-forward map
\begin{equation}\label{v123}
\cS_{g,n}^* \rightarrow A^*(\overline{\mathcal{M}}_{g,n}, \mathbb{Q})
\end{equation}
is a ring homomorphism,
 see \cite[Section 0.3]{PPZ13}
for a detailed discussion.

The image of \eqref{v123} is, by definition, the {\em tautological subring}
$$R^*(\overline{\mathcal{M}}_{g,n})\subset A^*(\overline{\mathcal{M}}_{g,n}, \mathbb{Q})\ .$$
Hence, we have a quotient
$$ \cS_{g,n}^*
\stackrel{q}{\longrightarrow} R^*(\overline{\mathcal{M}}_{g,n}) \longrightarrow 0\ .$$
The ideal of {\em tautological relations} is the kernel of $q$.

\subsection{Vertex, leg, and edge factors} 
Pixton's relations are determined by 
a set 
$$\tP=\{\cR_{g,A}^d\}$$ of
elements $\cR^d_{g,A}\in \cS^d_{g,n}$ associated to the data
\begin{enumerate}
\vspace{7pt}
\item[$\bullet$] $g,n\in \mathbb{Z}_{\geq 0}$ in the stable range $2g-2+n>0$,
\vspace{7pt}
\item[$\bullet$] $A=(a_1,\ldots, a_n), \ \ a_i \in\{0,1\}$,
\vspace{7pt}
\item[$\bullet$] $d\in \mathbb{Z}_{\geq 0}$ satisfying
$d > \frac{g-1+\sum_{i=1}^n a_i}{3}$.
\end{enumerate}
\vspace{7pt}
The elements $\cR^d_{g,A}$ are expressed as sums over
stable graphs of genus $g$ with $n$ legs. 
Before writing the formula for $\mathcal{R}^d_{g,A}$, a few definitions are required.

The form of the
hypergeometric series $\AAA$ and $\BBB$ used in
Pixton's relations is following:
\begin{align*}
\CCC_0(T)&=\AAA(-288T) = \sum_{i=0}^\infty \frac{(6i)!}{(2i)!(3i)!}(-T)^i
= 1-60T+27720T^2 -\cdots,\\
\CCC_1(T)&=-\BBB(-288T) = -\sum_{i= 0}^\infty  \frac{(6i)!}{(2i)!(3i)!}\frac{6i+1}{6i-1}
(-T)^i
=1 + 84T - 32760T^2 + \cdots.
\end{align*}
These series control the original Faber-Zagier relations
 and continue
to play a central role in the set $\tP$.

Let $f(T)$ be a power series with vanishing constant and linear terms,
$$f(T)\in T^2\mathbb{Q}[[T]]\ .$$
For each $\oM_{g,n}$,  we define
\begin{equation}\label{g33g}
\kappa(f) = \sum_{m \geq 0} \frac1{m!}\ { p_{m*}} \Big(f(\psi_{n+1}) \cdots f(\psi_{n+m})\Big)
 \ \in A^*(\oM_{g,n},\mathbb{Q}),
\end{equation}
where $p_m$ is the forgetful map
$$p_m: \oM_{g,n+m} \to \oM_{g,n}.$$ 
By the vanishing in degrees 0 and 1  of $f$, the sum \eqref{g33g} is finite.

Let $\mathsf{G}_{g,n}$ be the (finite) set of stable graphs of
genus $g$ with $n$ legs (up to isomorphism).
Let $\Gamma \in \mathsf{G}_{g,n}$. For each vertex $v\in \V$,
we introduce an auxiliary variable $\zeta_v$ and impose the
conditions
$$\zeta_v \zeta_{v'}= \zeta_{v'} \zeta_v\ , \ \ \ \zeta_v^2= 1\ .$$
 The variables $\zeta_v$ will be responsible for keeping track of
a local parity condition at each vertex.

The formula for $\mathcal{R}_{g,A}^d$ is a sum over $\mathsf{G}_{g,n}$.
The summand corresponding to $\Gamma \in \mathsf{G}_{g,n}$ is a 
product of 
vertex, leg, and edge factors:
\begin{enumerate}
\vspace{5pt}
\item[$\bullet$]
For $v\in \V$, let
$\kappa_v = \kappa\big(T-T \B_0(\zeta_vT)\big)$.
\vspace{5pt}
\item[$\bullet$]
For $l \in \L$,
let
 $\B_l =\zeta_{v(l)}^{a_l} \B_{a_l} \! \left(\zeta_{v(l)} \psi_{l}\right)$,  
where $v(l)\in V$ is the vertex to which the leg is assigned.
\vspace{0pt}
\vspace{-3pt}
\item[$\bullet$]
For $e\in \E$, let
\begin{align*}
\Delta_e &= \frac{ \zeta' + \zeta'' - 
\B_0(\zeta' \psi') \zeta''\B_1(\zeta'' \psi'')
-\zeta'\B_1(\zeta' \psi') \B_0(\zeta'' \psi'')}
{\psi'+\psi''}\\
&=(60 \zeta' \zeta''-84) +
\left[32760(\zeta'\psi' + \zeta'' \psi'') - 27720 (\zeta'\psi''+\zeta''\psi')\right] + \cdots,
\end{align*}
where $\zeta',\zeta''$ are the $\zeta$-variables assigned to the vertices adjacent to the edge $e$ and $\psi', \psi''$ are the $\psi$-classes corresponding to the half-edges.
\end{enumerate}
\vspace{5pt}
The numerator of $\Delta_e$ is divisible by its denominator due to the identity
$$
\B_0(T) \B_1(-T) + \B_0(-T) \B_1(T) =2.
$$
Certainly, $\Delta_e$ is symmetric in the
half-edges.

\subsection{Relations}

Let $A = (a_1, \dots, a_n) \in \{0,1\}^n$.
We denote by
$\cR_{g,A}^d\in \cS_{g,n}^d$ the degree $d$ component of the strata algebra class 
$$
\sum_{\Gamma\in \mathsf{G}_{g,n}} \frac1{|{\text{Aut}}(\Gamma)| }
\, 
\frac1{2^{h^1(\Gamma)}}
\;
\left[\Gamma, \; \Bigl[
\prod \kappa_v \prod \B_l 
\prod \Delta_e
\Bigr ]_{\prod_v \zeta_v^{\mathrm{g}(v)-1}}
\right] \ \in \cS_{g,n},
$$
where the products are taken over all vertices, all legs, and all edges of the graph~$\Gamma$.
The subscript $\prod_v \zeta_v^{\mathrm{g}(v)-1}$ indicates
the coefficient of the monomial $\prod_v \zeta_v^{\mathrm{g}(v)-1}$
after the product inside the brackets is expanded.

We denote by $\tP$ the set of classes $\cR^d_{g,d}$ for
$$
d > \frac{g-1 + \sum_{i=1}^n a_i}{3}.
$$
By the following result, $\tP$ is a set of tautological relations.

\begin{theorem}[Janda] \label{Thm:relations}
Every element $\cR_{g,A}^d\in \tP$ lies in the kernel of 
the homomorphism $$q:\cS_{g,n}^* \rightarrow A^*(\oM_{g,n},\mathbb{Q})\ .$$ 
\end{theorem}

Pixton's relations were conjectured first in \cite{Pix12} and
and first proven to hold in $H^*(\oM_{g,n},\mathbb{Q})$ in
 \cite{PPZ13} using Witten's 3-spin theory and the
Givental-Teleman classification of CohFTs. The 
Chow results of \cite{Jan13} are via a study of  the equivariant Gromov-Witten
theory of $\mathbb{CP}^1$. The application of
the virtual localization formula \cite{GrP} to the Gromov-Witten theory
of $\mathbb{CP}^1$ bypasses the need for Teleman's cohomological
results, and the proof of \cite{Jan13} holds in Chow.

Both Witten's 3-spin theory and the equivariant Gromov-Witten
theory of $\mathbb{CP}^1$ are {\em Chow Field Theories}: the CohFT
axioms are satisfied in the Chow rings $A^*(\oM_{g,n},\mathbb{Q})$.

\begin{question}
  Does the Givental-Teleman classification hold for semisimple Chow
  Field theories?
\end{question}

The proof of \cite{PPZ13} would imply the Chow vanishing of
Theorem \ref{Thm:relations} if the answer to Question 4 is yes.
The classification of CohFTs by Teleman uses the stable
cohomology of the moduli spaces of curves. Are such stability
results valid in the Chow ring?

\end{document}